\documentclass[11pt]{article}
\usepackage{amssymb}
\usepackage{amsmath}
\usepackage{amsthm}
\newcommand{\bnu}{\boldmath$\nu$}
\newcommand{\mbnu}{\mbox{\bnu}}
\newcommand{\bomega}{\boldmath$\omega$}
\newcommand{\mbomega}{\mbox{\bomega}}
\newcommand{\bDelta}{\boldmath$\Delta$}
\newcommand{\mbDelta}{\mbox{\bDelta}}

\newcommand{\bGamma}{\boldmath$\Gamma$}
\newcommand{\mbGamma}{\mbox{\bGamma}}
\newcommand{\bSigma}{\boldmath$\Sigma$}
\newcommand{\mbSigma}{\mbox{\bSigma}}
\newcommand{\btheta}{\boldmath$\theta$}
\newcommand{\mbtheta}{\mbox{\btheta}}

\newcommand{\reg}{$\emptyset \subseteq a \subseteq P$}
\newcommand{\rag}{$P_{1} \subseteq a \subseteq P$}
\newcommand{\Xbar}{$\bar{\bf X}$}
\newcommand{\Ybar}{$\bar{\bf Y}$}
\newcommand{\mXbar}{\mbox{\Xbar}}
\newcommand{\mYbar}{\mbox{\Ybar}}
\newcommand{\Zbar}{$\bar{\bf Z}$}
\newcommand{\mZbar}{\mbox{\Zbar}}

\newfont{\bfit}{cmbxti10 scaled\magstep1}

\newcommand{\mbS}{{\bf S}}

\newcommand{\mbg}{\boldsymbol{g}}

\DeclareMathOperator{\tr}{tr}

\begin{document}

\newtheorem{theorem}{\indent {\sc Theorem}}
\newtheorem{prop}{\indent {\sc Proposition}}
\newtheorem{lemma}{\indent {\sc Lemma}}
\renewcommand{\proofname}{\hspace*{\parindent}{\sc Proof.}}
\thispagestyle{empty}

\begin{center}
{\bf Towards a unified theory for testing statistical hypothesis: Multinormal mean with nuisance
covariance matrix}\\
\medskip
{Ming-Tien Tsai} \\
\medskip
{\it Academia Sinica, Taipei, Taiwan}
\end{center}

\vspace{0.3cm}
\noindent {\bf Abstract}
\vspace{0.2cm}

\begin{center}
\small
\parbox{15cm}{\sloppy \ \,
Under a multinormal distribution with an arbitrary unknown covariance matrix, the main purpose of this paper is to propose
a framework to achieve the goal of reconciliation of Bayesian, frequentist, and Fisher's reporting $p$-values,
Neyman-Pearson's optimal theory and Wald's decision theory for the problems of testing mean against restricted alternatives
(closed convex cones). To proceed, the tests constructed via the likelihood ratio (LR) and the union-intersection (UI)
principles are studied. For the problem of testing against the positive orthant space alternative, neither Fisher's approach
of reporting $p$-values alone nor Neyman-Pearson's optimal theory for power function alone is a satisfactory criterion for
evaluating the performance of tests. We show that the LRT and the UIT are not the proper Bayes tests, however, we also
show that the LRT and the UIT are the integrated LRT and the integrated UIT, respectively. Wald's decision theory via
$d$-admissibility may shed light on resolving these challenging issues of imposing the balance between type 1 error and
power. We generalize the method of convex acceptance region of Birnbaum-Stein to show that for the problems of testing
against restricted alternatives, the UITs are $d$-admissible. The hyperbolic type set of acceptance regions of LRTs may
cause difficulty in studying whether the LRTs are $d$-admissible.}
\end{center}

\vspace{.3cm}
\indent{\small{\it AMS 2001 subject classification}: 62C15, 62H15.}


\vspace{0.3cm}
\noindent {\bf 1. Introduction and preliminary notions}
\vspace{0.2cm}

\def \theequation{1.\arabic{equation}}
\setcounter{equation}{0}

\indent Let \{${\bf X}_{i}; 1 \leq i \leq n$\} be independent and identically distributed random vectors (i.i.d.r.v.)
having a $p$-variate normal distribution with mean vector ${\mbtheta}$ and unknown dispersion matrix ${\mbSigma} \in
{\cal M}$, which denotes the set of all the $p \times p$ positive definite (p.d.) matrices. Consider the hypotheses
\begin{eqnarray}
H_{0}:{\mbtheta} = {\bf 0}~~\hbox{vs.}~~H_{1}:{\mbtheta} \in {\cal C} \backslash \{{\bf 0}\},
\end{eqnarray}
where ${\cal C}$ denotes the set of closed convex cones for which every convex cone containing a p-dimensional open set.
Denote the positive orthant space by ${\cal O}_{p}^{+} = \{{\mbtheta}\in R^{p}|~{\mbtheta} \geq {\bf 0}\}$ and the
half-space by ${\cal H}_{p}^{*}=$

\vspace{0.3cm}
\noindent \underline{~~~~~~~~~~~~~~~~~~~~~~~~~~~~~~~~~~~~~~~~~~~~}\\
E-mail: mttsai@stat.sinica.edu.tw
\pagebreak

\noindent $\{{\mbtheta}\in R^{p}|~{\theta}_{p} \geq 0 \}$ throughout this paper. Let ${\cal S}_{1}$
be the set of half-spaces, and ${\cal S}={\cal C}\backslash {\cal S}_1$. Notice that when a closed convex cone within
${\cal S}$ is a proper set contained within a half-space, under a suitable linear transformation the problem in (1.1) can
then be reduced to the problem for testing against the positive orthant space alternative with another unknown positive
definite covariance matrix and the transformed data set. When ${\cal C}$ is a specific half-space, it can be transformed
into another half-space using a non-singular linear transformation based on the transformed data. Therefore, for simplicity
it is sufficient to examine scenarios where ${\cal C}$ represents the positive orthant space ${\cal O}_{p}^{+}$ and where
${\cal C}$ is the half-space ${\cal H}_{p}^{*}$.

\indent The likelihood ratio test (LRT) is constructed using the likelihood ratio principle, while the union-intersection
test (UIT) is constructed based on the union-intersection principle (Roy [21]). For testing against the global alternative,
$H^{g}_{0}:{\mbtheta} = {\bf 0}~~\hbox{vs.}~~H^{g}_{1}:{\mbtheta} \ne {\bf 0}$, the LRT and the UIT are isomorphic, and they
are well-known as Hotelling's $T^{2}$-test (Anderson [1]).

\indent For testing against the global alternative, Hotelling's $T^{2}$-test enjoys many optimal properties of the
Neyman-Pearson testing hypothesis theories such as similarity, unbiasedness, power monotonicity, most stringent, uniformly
most powerful invariant and $\alpha$-admissibility (Anderson [1]). However, when addressing the problem of testing against
the positive orthant space alternative, both the LRT and the UIT are demonstrated to be biased (Perlman [18], Sen and Tsai
[24]). Furthermore, it is shown that both the LRT and the UIT are power-dominated by the corresponding LRT and UIT for the
problem of testing against the half-space alternative ${\cal H}^{*}_{p}$, respectively. As such, both the LRT and the UIT
for the problem of testing against the positive orthant space alternative are $\alpha$-inadmissible. The related power
domination problems are also studied in Section 2, some power domination results are against the common statistical sense.
One of the primary reasons for this is that for the problem of testing against the positive orthant space alternative, the
distribution functions of both the LRT statistic and the UIT statistic are under the dependence on the unknown nuisance
covariance matrix under the null hypothesis. Because of this unpleasant feature, the Fisher's approach to reporting
$p$-values may also not work well for the problem of testing against the positive orthant space alternative. Several
methods have been proposed to eliminate the unpleasant features, we discussed some in Section 2.

\indent For testing against the global alternative, Kiefer and Schwartz [11] demonstrated that Hotelling's $T^{2}$-test is
a proper Bayes test. Additionally, Hotelling's $T^{2}$-test is also known as a version of the integrated LRT for objective
priors (Berger, Liseo, and Wolpert [3]). There are several ways to establish the equivalent relationships of Hotelling's
$T^{2}$-test with Bayes tests. On the other hand, for the problems of testing against restricted alternatives considered
in (1.1), the LRT and the UIT differ (Perlman [18], Sen and Tsai [24]) when the covariance matrix ${\mbSigma}$ is
unknown. The problem (1.1) provides us some partial information about the parameters, the difficulty of integration over
a high-dimensional restricted parameter space, the explicit forms of Bayes tests are hardly derived though the numerical
values of them can be obtained by the Markov Chain Monte Carlo method. In Section 3, we will demonstrate that the
LRT and the UIT are no longer to be the proper Bayes tests for the problem (1.1), while it is shown that both the LRT
and the UIT are the versions of the integrated LRT ans the integrated UIT for problem (1.1), respectively. This leads to
the reconciliation of the frequentist and Bayesian paradigms. Hence, it provides us a hope to further investigate
the optimal propertis of UIT and LRT.

\indent Elimination of nuisance parameters is a central issue in statistical inference, there are many ways to eliminate
the nuisance parameters in the literature (Basu [2]). Most of these methods can be applied to the problem of testing against
the positive orthant space alternative. In Section 4, we show that Hotelling's $T^{2}$-tests are inadmissible for the
problem (1.1). Furthermore, we also show that the UIT is $d$-admissible for testing against the positive orthant space
alternative, and that it is $\alpha$-admissible for testing against the half-space alternative. However, we fail to
claim these results for LRTs due to the facts that the acceptance regions of LRTs are the hyperbolic type sets.

\vspace{0.3cm}
\def \theequation{2.\arabic{equation}}
\setcounter{equation}{0}

\vspace{0.3cm}
\noindent {\bf 2. Power comparision and/or reporting $p$-values?}
\vspace{0.2cm}

\indent For every $n$ $(\geq 2)$, let
\begin{eqnarray}
  {\mXbar}_{n} = n^{-1} \sum_{i=1}^{n} {\bf X}_{i}~~\hbox{and}~~{\bf S}_{n} =(n-1)^{-1}
  \sum_{i=1}^{n} ({\bf X}_{i} - {\mXbar}_{n}) ({\bf X}_{i} -
  {\mXbar}_{n})',
\end{eqnarray}
the Hotelling's $T^{2}$-test is expressed as
\begin{eqnarray}
T^{2}=n{\mXbar}'_{n}{\bf S}^{-1}_{n}{\mXbar}_{n}.
\end{eqnarray}

\indent Historically, two different approaches are adopted in hypothesis testing problems to evaluate the tests for
the frequientist paradigm: one is Fisher's approach to report $p$-values of the test, and the other is Neyman-Pearson's
approach being used to determine the power of test at a fixed significance level. Both the Fisher approach and the
Neyman-Pearson approach are widely adapted to measure the performance of Hotelling's $T^{2}$-test for the problem
of testing against the unrestricted alternative.

\vspace{0.3cm}
\noindent {\it 2.1 The LRT and the UIT}
\vspace{0.2cm}

\indent Let $P =\{1,\ldots, p\}$, and for every $a$: {\reg}, let $a'$ be its complement and $|a|$ its
cardinality. For each $a$, we partition ${\mXbar}_{n}$ and ${\bf S}_{n}$ as
\begin{eqnarray}
  {\mXbar}_{n}=\left(
  \begin{array}{c}
    {\mXbar}_{na}\\
    {\mXbar}_{na'}
  \end{array}
  \right) ~~\hbox{and}~~
  {\bf S}_{n} = \left(
  \begin{array}{cc}
    {\bf S}_{naa}& {\bf S}_{naa'}\\
    {\bf S}_{na'a}& {\bf S}_{na'a'}
  \end{array}
  \right),
\end{eqnarray}
and write
\begin{eqnarray}
  {\mXbar}_{na:a'} = {\mXbar}_{na}-{\bf S}_{naa'}{\bf S}_{na'a'}^{-1}
  {\mXbar}_{na'}, \\
  {\bf S}_{naa:a'}={\bf S}_{naa}-{\bf S}_{naa'}{\bf S}_{na'a'}^{-1}
  {\bf S}_{na'a}.
\end{eqnarray}
Further, let
\begin{eqnarray}
  I_{na} = 1\{{\mXbar}_{na:a'}>{\bf 0}, {\bf S}_{na'a'}^{-1}{\mXbar}_{na'}\leq 0\},
\end{eqnarray}
for {\reg}, where $1\{\cdot\}$ denotes the indicator function.
Then for the problem of testing $H_{0}: {\mbtheta} = {\bf 0}$ against the positive orthant
space alternative $H_{1\it{O}^{+}}: {\mbtheta} \in {\cal O}^{+}_{p} \backslash \{{\bf 0}\}$, from the
results of Perlman [18] and Sen and Tsai [24] the LRT and the UIT statistics are of the
forms
\begin{eqnarray}
  L_{n} = \sum_{\mbox{\reg}}\left\{ {\frac{n{\mXbar}'_{na:a'}
  {\bf S}_{naa:a'}^{-1} {\mXbar}_{na:a'}} {1+n{\mXbar}'_{na'} {\bf S}_{na'a'}^
  {-1} {\mXbar}_{na'}}} \right\}I_{na}
\end{eqnarray}
and
\begin{eqnarray}
  U_{n}= \sum_{\mbox{\reg}}\{n{\mXbar}'_{na:a'}{\bf S}_{naa:a'}^
 {-1}{\mXbar}_{na:a'}\}I_{na}
\end{eqnarray}
respectively.

\indent Side by side, for testing $H_{0}: {\mbtheta} = {\bf 0}$ against the half-space alternatives
$H_{1\it{H}^{*}}: {\mbtheta} \in {\cal H}^{*}_{p}\backslash \{{\bf 0}\}$ the LR and the UI
test statistics are of the forms
\begin{eqnarray}
  L_{n}^{\ast} = \sum_{\mbox{\rag}}\left\{ {\frac{n{\mXbar}'_{na:a'}
  {\bf S}_{naa:a'}^{-1} {\mXbar}_{na:a'}} {1+n{\mXbar}'_{na'} {\bf S}_{na'a'}^
  {-1} {\mXbar}_{na'}}} \right\}I^{*}_{na}
\end{eqnarray}
and
\begin{eqnarray}
  U_{n}^{\ast}= \sum_{\mbox{\rag}}\{n{\mXbar}'_{na:a'}{\bf S}_{naa:a'}^
 {-1}{\mXbar}_{na:a'}\}I^{*}_{na}
\end{eqnarray}
respectively, where $P_{1}=P-\{ p\}, I^{*}_{nP}=1\{\Bar{X}_{np}>0\}$ and
$I^{*}_{nP1}=1\{\Bar{X}_{np}\leq 0\}$.

\indent It is clear from (2.7)-(2.10) that
\begin{eqnarray}
 L^{*}_{n} \geq L_{n}~~\hbox{and}~~U^{*}_{n} \geq U_{n}
\end{eqnarray}
with probability one.

\indent We would like to determine whether the Fisher approach and/or the Neyman-Pearson approach can be well adapted
for the problem (1.1). First, we consider the problem of testing against the positive orthant space alternative.

\vspace{0.3cm}
\noindent {\it 2.2 Testing against the positive orthant space alternative ${\cal O}_{p}^{+}$}
\vspace{0.2cm}

\indent  For the problem of testing against the positive orthant space alternative, it follows from the arguments of
Perlman [18] and Sen and Tsai [24] that for every c $>$ 0,
\begin{eqnarray}
  P\{ L_{n} \geq c| H_0,{\mbSigma} \} = P_{{\bf 0},{\mbSigma}} \{ L_{n} \geq c \} \\
  = \sum_{k=1}^{p}w(p,k;{\mbSigma})P\{ \chi_{k}^{2} /\chi_{n-p}^{2}\geq c\}\nonumber
\end{eqnarray}
and
\begin{eqnarray}
  P\{ U_{n} \geq c| H_0,{\mbSigma} \} = P_{{\bf 0},{\mbSigma}} \{ U_{n} \geq c \} \\
  = \sum_{k=1}^{p}w(p,k;{\mbSigma}) {\bar G}_{n,k,p}^{\ast}(c),\nonumber
\end{eqnarray}
respectively, where for each $k$ (=0,1,\ldots,$p$)
\begin{eqnarray}
  w(p,k;{\mbSigma})=\sum\nolimits_{\{a\subseteq P: |a|=k\}}P\{{\bf Z}_{a:a'}>{\bf 0},
  {\mbSigma}_{a'a'}^{-1}{\bf Z}_{a'} \leq {\bf 0}\},
\end{eqnarray}
and ${\bf Z}\sim{\it N}_{p}({\bf 0},{\mbSigma})$ and the partitioning is made
as in (2.3)-(2.5) (with ${\bf S}_{n}$ replaced by {\bSigma}).

\indent By virtue of (2.12), (2.13) and (2.14), we encounter the problem for which both the null distributions of LRT and UIT
statistics depend on the unknown covariance matrix, though the dependence of (2.12) and (2.13) on {\bSigma} is only
through the $w(p,k;{\mbSigma})$. Hence, we have difficulty adopting Fisher's approach to reporting the $p$-values. On the
other hand, according to the definition of level of significance, Perlman [18] and Sen and Tsai [24] based on the
Neyman-Pearson optimal property allowing ${\mbSigma}$ to vary over the entire class ${\mathcal M}$ and obtained that for
every c $> 0$,
\begin{eqnarray}
  \sup_{\{{\mbSigma}\in {\mathcal M}\}}P_{{\bf 0},{\mbSigma}}\{L_{n}\geq c\}
  = {\frac {1}{2}}\left[\bar G_{p-1,n-p}(c)+\bar G_{p,n-p}(c)\right]
\end{eqnarray}
and
\begin{eqnarray}
  \mathop{\sup}_{\{{\mbSigma}\in {\mathcal M}\}}P_{{\bf 0},{\mbSigma}}\{U_{n} \geq
  c\}={\frac{1}{2}}\left [\bar G_{n,p-1,p}^{*}(c)+\bar G_{n,p,p}^{*}(c)\right],
\end{eqnarray}
respectively.

\indent The right hand side of (2.15) and (2.16) are equated to $\alpha$, we then can get the critical values of LRT and UIT,
respectively.

\vspace{0.3cm}
\noindent {\it 2.3 Testing against the half-space alternative ${\cal H}_{p}^{*}$}
\vspace{0.2cm}

\indent Let $\chi_{m}^{2}$ be the central chi-square distribution with $m(\geq 0)$ degrees of
freedom. Denoted by
\begin{eqnarray}
  G_{a,b}(u)=P\{\chi_{a}^{2}/\chi_{b}^{2}\leq u\}, ~u\in R^{+},
\end{eqnarray}
and ${\bar G}_{a,b}(u)=1-G_{a,b}(u)$, where $R^{+}$ denotes the positive real number.  Also
let
\begin{eqnarray}
  {\bar G}_{n,a,p}^{\ast}(u) = \int_{0}^{\infty} {\bar G}_{a,n-p}({\frac{u}{1+t
  }})dG_{p-a,n-p+a}(t),~~u\in R^{+},
\end{eqnarray}
which is the convolution of the d.f.'s of $\chi_{a}^{2}/\chi_{n-p}^{2}$ and 1+$\chi_{p-a}^{2}/\chi_{n-p+a}^{2}$.
By (2.15) and (2.16), as such for the problem of testing against the half-space alternative ${\mathcal H}^{*}_{p}
\backslash \{{\bf 0}\}$, it is easy to show that when ${\mbSigma} \in {\cal M},~\forall ~c > 0$
\begin{eqnarray}
 P_{{\bf 0}, {\mbSigma}}\{L_{n}^{\ast}\geq c\}
 &=& {\frac{1}{2}}\left [\bar G_{p-1,n-p}(c)+\bar G_{p,n-p}(c)\right],
\end{eqnarray}
and
\begin{eqnarray}
 P_{{\bf 0}, {\mbSigma}}\{U_{n}^{\ast}\geq c\}
  &=&{\frac{1}{2}}\left [\bar G_{n,p-1,p}^{*}(c)+\bar G_{n,p,p}^{*}(c)\right].
\end{eqnarray}
Notice that (2.19) and (2.20) are free from the nuisance parameter ${\mbSigma}$, namely, the $p$-values of LRT and UIT do
not depend on the nuisance parameter. By the result of Tang [28] and similar arguments in Theorems 2.1 and 3.1 of
Sen and Tsai [24] we then obtain that both LRT and UIT for testing against the half-space alternative are similar and
unbiased, respectively. Hence, for the problem of testing against the half-space alternative both approaches of Fisher and
Neyman-Pearson can well explain the performance of the LRT and the UIT.
\vspace{0.2cm}

\vspace{0.3cm}
\noindent {\it 2.4 Power domination}
\vspace{0.2cm}

\indent We introduce some notations first, let $\pi_{\bf A}({\bf x};{\cal C})$ be the orthogonal projection of ${\bf x}$
onto ${\cal C}$ with respect to the inner product $<, >_{\bf A}$, then
\begin{eqnarray}
   U_{n}=|\!|{\pi}_{{\bf S}_n}(n^{1/2}{\mXbar}_{n}; {\cal O}_{p}^{+})|\!|_{{\bf S}_n}^{2}
\end{eqnarray}
and
\begin{eqnarray}
   U^{*}_{n}=|\!|{\pi}_{{\bf S}_n}(n^{1/2}{\mXbar}_{n}; {\mathcal H}^{*}_{p})|\!|_{{\bf S}_n}^{2}.
\end{eqnarray}

\indent And the LRT statistics in (2.7) and (2.9) can then be represented as
\begin{eqnarray}
  L_{n}= |\!| {\pi}_{{\bf S}_{n}}(n^{1/2} {\mXbar}_{n}; {\cal O}^{+}_{p})
    |\!|_{{\bf S}_{n}}^{2} \{ 1 + |\!| n^{1/2} {\mXbar}_{n} - {\pi}_{{\bf S}_{n}}
     (n^{1/2} {\mXbar}_{n}; {\cal O}^{+}_{p} ) |\!|_{{\bf S}_{n}}^{2} \}^{-1}
\end{eqnarray}
and
\begin{eqnarray}
   L^{*}_{n}=|\!| {\pi}_{{\bf S}_{n}}(n^{1/2} {\mXbar}_{n}; {\cal H}^{*}_{p})
    |\!|_{{\bf S}_{n}}^{2} \{ 1 + |\!| n^{1/2} {\mXbar}_{n} - {\pi}_{{\bf S}_{n}}
     (n^{1/2} {\mXbar}_{n}; {\cal H}^{*}_{p} ) |\!|_{{\bf S}_{n}}^{2} \}^{-1}
\end{eqnarray}
respectively.

\indent Based on the results of (2.15) and (2.19), we will easily find that both LRTs for the problems of testing against
the positive orthant space altenative and of testing against the half-space alternative have the same critical value when
the maximum principle is used to define the level of significance $\alpha$. So do the results of (2.16) and (2.20) for UITs.
These phenomena will bring the results which are against our common statistical sense. First, we have

\vspace{0.3cm}
\noindent {\bf Proposition 2.1.} {\it The UIT for testing } $H_{0}:{\mbtheta}={\bf 0}$ against
$H_{1}: {\mbtheta} \in {\cal C}\backslash \{{\bf 0}\}$, {\it where} ${\cal C}$ {\it is a closed convex cone such that}
${\cal U} \subseteq {\cal C} \subseteq {\mathcal H}^{*}_{p}$ {\it with } ${\cal U}$ {\it being a p-dimensional open set
and being contained in the half-space ${\mathcal H}^{*}_{p}$, {\it is power-dominated by the UIT for testing}
$H_{0}:{\mbtheta} ={\bf 0}~{against}~H_{1}: {\mbtheta} \in {\mathcal H}^{*}_{p} \backslash \{{\bf 0}\}$.}

\noindent {\bf Proof}. First note that ${\mathcal C} \subseteq {\mathcal H}^{*}_{p}$, hence
\begin{eqnarray}
 |\!| \pi_{{\bf S}_n}(n^{1/2}{\mXbar}_{n}; {\mathcal H}^{*}_{p}) |\!|_{{\bf S}_n}^{2}
  \geq |\!| \pi_{{\bf S}_n}(n^{1/2}{\mXbar}_{n}; {\cal C}) |\!|_{{\bf S}_n}^{2}. \nonumber
\end{eqnarray}
Then for every $c >0,$
\begin{eqnarray}
   P_{{\mbtheta},{\mbSigma}} \{|\!|{\pi}_{{\bf S}_n}(n^{1/2}{\mXbar}_{n};
   {\mathcal H}^{*}_{p})|\!|_{{\bf S}_n}^{2}  \geq c \} \geq P_{{\mbtheta},{\mbSigma}}
   \{ |\!| \pi_{{\bf S}_{n}}(n^{1/2}{\mXbar}_{n}; {\mathcal C}) |\!|_{{\bf S}_{n}}^{2}
   \geq c \}. \nonumber
\end{eqnarray}
\indent Let ${\mbSigma}^{-1}
={\bf B}^{'}{\bf B}$, and ${\mYbar}_{n}={\bf B}{\mXbar}_{n}$, and ${\bf T}_{n}
={\bf B}{\bf S}_{n}{\bf B}^{'}.$ Since ${\bf B}{\mathcal H}^{*}_{p}={\mathcal H}^{+}_{p}$
(another half-space), thus $P_{{\bf 0},{\mbSigma}} \{|\!|{\pi}_{{\bf S}_n}(n^{1/2}
{\mXbar}_{n}; {\mathcal H}_{p}^{+})|\!|_{{\bf S}_n}^{2} \geq c \} = P_{{\bf 0},{\bf I}}
\{|\!|{\pi}_{{\bf T}_n}(n^{1/2}{\mYbar}_{n};{\mathcal H}_{p}^{+})|\!|_{{\bf T}_n
  }^{2} \geq c \}$. Notice that ${\cal E}({\bf Y}_{n})=$(n-1)${\bf I}$,
thus $P_{{\bf 0},{\mbSigma}} \{|\!|{\pi}_{{\bf S}_n}(n^{1/2} {\mXbar}_{n};
{\mathcal H}_{p}^{+})|\!|_{{\bf S}_n}^{2} \geq c \}$ is independent of the unknown
${\mbSigma}$. Consider the sequence $\{{\mbSigma}_{k}\}$ and let ${\mbSigma}_{k}^{-1}
={\bf B}_{k}^{'}{\bf B}_{k}$, for all $k \geq 1$, such that ${\bf B}_{k}
({\mathcal C_{\lambda}}) \subseteq {\bf B}_{k+1}({\cal C_{\lambda}})$ and
$\cup_{k=1}^{\infty}{\bf B}_{k} ({\mathcal C_{\lambda}}) =int({\cal H}_{p}^{+})$, where
${\mathcal C_{\lambda}}$ with $0 <{\lambda} < 1$ is a right circular cone defined as in
$2.3^{o}$ of Perlman [18]. The same proof as in Theorem 6.2 of Perlman [18], we obtain
\begin{eqnarray}
 & &P_{{\bf 0},{\mbSigma}_{k}} \{|\!|{\pi}_{{\bf S}_n}(n^{1/2}{\mXbar}_{n};
{\mathcal H}_{p}^{+})|\!|_{{\bf S}_n }^{2} \geq c \} \geq
 P_{{\bf 0},{\mbSigma}_{k}}\{ |\!| \pi_{{\bf S}_{n}} (n^{1/2}{\mXbar}_{n};
  {\mathcal C}) |\!|_{{\bf S}_{n}}^{2} \geq c \}\nonumber \\
   & &\geq  P_{{\bf 0},{\mbSigma}_{k}}\{ |\!| \pi_{{\bf S}_{n}} (n^{1/2}{\mXbar}_{n};
  {\mathcal C_{\lambda}}) |\!|_{{\bf S}_{n}}^{2} \geq c \} =
  P_{{\bf 0},{\bf I}}\{ |\!| \pi_{{\bf S}_{n}^{*}} (n^{1/2}{\mXbar}_{n}^{*};
  {\bf B}_{k}({\mathcal C_{\lambda}})) |\!|_{{\bf S}_{n}^{*}}^{2} \geq c \},\nonumber
\end{eqnarray}
where ${\mXbar}_{n}^{*}={\bf B}_{k}{\mXbar}_{n}$, and ${\bf S}^{*}_{n}=
{\bf B}_{k}{\bf S}_{n}{\bf B}_{k}^{'}$. By Lemma 2.8 of Brown [7], we may conclude that
the power functions of the UIT for testing $H_{0}:{\mbtheta}={\bf 0}~\hbox{vs.}~H_{1}:
{\mbtheta} \in {\mathcal C}$ are analytic on the space $\{{\mbtheta} \in {\mathcal C},
{\mbSigma} > {\bf 0}\}$, and hence they are continuous. Thus
$\lim_{k \to \infty}P_{{\bf 0},{\bf I}}\{ |\!| \pi_{{\bf S}^{*}_{n}}
   (n^{1/2}{\mXbar}_{n}^{*}; {\bf B}_{k}({\cal C_{\lambda}})) |\!|_{{\bf S}^{*}_{n}}^{2}
   \geq c \}= P_{{\bf 0},{\bf I}}\{ |\!| \pi_{{\bf S}_{n}^{*}} (n^{1/2}{\mXbar}_{n}^{*};
    {\cal H}_{p}^{+}) |\!|_{{\bf S}^{*}_{n}}^{2} \geq c \}$, and hence
\begin{eqnarray}
{\sup}_{\{{\mbSigma} \in
{\cal M}\}} P_{{\bf 0},{\mbSigma}}\{ |\!| \pi_{{\bf S}_{n}} (n^{1/2}{\mXbar}_{n};
  {\cal C}) |\!|_{{\bf S}_{n}}^{2} \geq c \}
  =P_{{\bf 0},{\mbSigma}} \{|\!|{\pi}_{{\bf S}_n}(n^{1/2}{\mXbar}_{n};
    {\cal H}_{p}^{+})|\!|_{{\bf S}_n }^{2} \geq c \},
\end{eqnarray}
i.e., the equation (2.16) is equivalent to the equation (2.20). Both UITs $U_{n}$ and $U^{*}_{n}$ have the same critical
value. Thus, the proposition follows.

\indent As for the LRT, note that for any ${\mathcal C} \subseteq {\mathcal H}^{*}_{p}$, $ |\!| \pi_{{\bf S}_n}(n^{1/2}
 {\mXbar}_{n}; {\mathcal H}^{*}_{p}) |\!|_{{\bf S}_n}^{2} \geq |\!| \pi_{{\bf S}_n}(n^{1/2}{\mXbar}_{n}; {\cal C})
  |\!|_{{\bf S}_n}^{2}$, and hence $|\!| n^{1/2} {\mXbar}_{n} - {\pi}_{{\bf S}_{n}}(n^{1/2} {\mXbar}_{n}; {\cal H}^{*}_{p} )
|\!|_{{\bf S}_{n}}^{2} \leq |\!| n^{1/2} {\mXbar}_{n} - {\pi}_{{\bf S}_{n}}(n^{1/2} {\mXbar}_{n}; {\cal C})
|\!|_{{\bf S}_{n}}^{2}$. Thus, by similar argument as in the proof of Proposition 2.1, we also have the following.

\vspace{0.3cm}
\noindent {\bf Proposition 2.2.} {\it The LRT for testing } $H_{0}:{\mbtheta}={\bf 0}~against
~H_{1}: {\mbtheta} \in {\cal C}\backslash \{{\bf 0}\}$, {\it where} ${\cal C}$ {\it is a closed convex cone such that}
${\cal U} \subseteq {\cal C} \subseteq {\mathcal H}^{*}_{p}$ {\it with } ${\cal U}$ {\it being a p-dimensional open set
and being contained in the half-space} ${\mathcal H}^{*}_{p}$, {\it is power-dominated by the LRT for testing}
$H_{0}:{\mbtheta} ={\bf 0}~{against}~H_{1}: {\mbtheta} \in {\mathcal H}^{*}_{p}\backslash \{{\bf 0}\}$.
\vspace{0.2cm}

\indent When using the maximum principle to define the level of significance $\alpha$, Proposition 2.2 indicates the
phenomenon, which was first noted by Tang [28], that the LRT for testing against the positive orthant space alternative is
power-dominated by that of testing against the half-space alternative. The result is against our common statistical sense.

\vspace{0.3cm}
\noindent {\it 2.5 Other ways to define the level of significance $\alpha$}
\vspace{0.2cm}

\indent To overcome these unpleasant phenomena mentioned above, instead of using the conservative maximum principle to
define the level of significance $\alpha$, i.e., taking the supremum of rejected probability under null hypothesis over
the set ${\cal M}$, we may use the Bayesian notions in this subsection, by taking the average of rejected probability
under the null hypothesis over the set ${\cal M}$ with respect to the weight function $g({\mbSigma})$ of ${\mbSigma}$,
to define the level of significance $\alpha$.

\indent The joint density function of $({\mXbar}_{n},{\bf S}_{n})$ under the null hypothesis is
\begin{eqnarray}
f({\mXbar}_{n}, {\bf S}_{n})=k_{0}|{\mbSigma}|^{-\frac{1}{2}n}
         |{\bf S}_{n}|^{\frac{1}{2}(n-p-2)}
         \mbox{e}^{{-\frac{1}{2}}\mbox{tr}({\mbSigma}^{-1}{\bf V}_{n})},
\end{eqnarray}
where $k_{0}$ is the normalizing constant and
${\bf V}_{n}=(n-1){\bf S}_{n}+n{\mXbar}_{n}{\mXbar}^{'}_{n}$.
First consider the inverted Wishart distribution $W^{-1}({\mbGamma},m)$ (Anderson [1]),
which is a proper prior of ${\mbSigma}$, as the weight function. Then we have
\begin{eqnarray}
 h_{1}({\mXbar}_{n},{\bf S}_{n})&=&\int_{{\mbSigma}\in {\mathcal M}}
     f({\mXbar}_{n}, {\bf S}_{n})g({\mbSigma}) ~d{\mbSigma} \\ \nonumber
 &=&  k_{2}~|{\mbGamma}|^{\frac{1}{2}m}|{\bf S}_{n}|^{\frac{1}{2}(n-p-2)}
        |{\bf V}_{n}+{\mbGamma}|^{-\frac{1}{2}(n+m-1)},
\end{eqnarray}
where the constant $k_{2}$ can be determined by the equation
$\int_{{\cal D}} ~h_{1}({\mXbar}_{n},{\bf S}_{n})~~dx_{1}\cdots dx_{p}=1$ with
${\cal D}=\{(x_{1},\cdots,x_{p})|~{\bf S}_{n}\in {\mathcal M} ~\hbox{in probability}\}$.
Let
\begin{eqnarray}
  b_{1}(k,n,p)=\sum_{|a|=k}\int_{I_{na}\cap {\cal D}}~h_{1}({\mXbar}_{n},{\bf S}_{n})
   ~dx_{1}\cdots dx_{p},
\end{eqnarray}
then it is easy to note that $ b_{1}(k,n,p) > 0, \forall k=1, 2, \cdots,p$ and
$\sum_{k=0}^{p}b_{1}(k,n,p)=1$.
Thus by Fubini theorem the critical values of the  LRT and the UIT for testing against the
positive orthant space alternative can then be determined via the formula
\begin{eqnarray}
  \alpha=\sum_{k=0}^{p}b_{1}(k,n,p){\bar G}_{k,n-p}(d^{1}_{n,{\alpha}})
\end{eqnarray}
and
\begin{eqnarray}
  \alpha=\sum_{k=0}^{p}b_{1}(k,n,p){\bar G}_{n,k,p}^{\ast}(c^{1}_{n,{\alpha}})
\end{eqnarray}
respectively. If the critical values of LRT and UIT are determined by the equations (2.29) and
(2.30), then obviously both the LRT and the UIT are similar and unbiased. These conclusions
are contrary to the ones made by Perlman [18] and Sen and Tsai [24] that both the LRT and
the UIT are not similar and biased.

\indent The level of significance $\alpha$ may be defined by using the noninformative priors
such as the Haar measure or the Bernardo reference prior, in which the corresponding posterior
densities exist, as the weight functions. Suppose that the Haar measure of ${\mbSigma}$ is
taken as the weight function, then we have
\begin{eqnarray}
  h_{2}({\mXbar}_{n},{\bf S}_{n})&=&k_{0}\int_{{\mbSigma}\in {\mathcal M}}~
     f({\mXbar}_{n}, {\bf S}_{n})|{\mbSigma}|^{-\frac{1}{2}(p+1)} ~~d{\mbSigma}   \\
   &=&k_{3}{|{\bf S}_{n}|}^{\frac{1}{2}(n-p-2)}{|{\bf V}_{n}|}^{-\frac{1}{2}n},  \nonumber
\end{eqnarray}
where the constant $k_{3}$ can be determined by the equation
$\int_{{\cal D}} ~h_{2}({\mXbar}_{n},{\bf S}_{n})~~dx_{1}\cdots dx_{p}=1$.
Let $b_{2}(k,n,p)$ be defined the same as in (2.12) with $h_{2}({\mXbar}_{n},{\bf S}_{n})$
instead of $h_{1}({\mXbar}_{n},{\bf S}_{n})$. As such, we may obtain another critical points
for the LRT and the UIT. Further note that condition on ${\bf V}_{n}$, then (2.31) reduces
to that of Wang and McDermott [32] for finding the critical value of semi-conditional LRT.

\indent From the above arguments, we can conclude that the phenomena regarding the LRT and the UIT can vary depending on
the ways of defining the level of significance $\alpha$. In passing, we note that vary weight functions yield different
critical values varying. The powers of the LRT and the UIT, respectively, heavily depend on the choice of weight functions
for the nuisance parameter ${\mbSigma}$. Therefore, for a specific test one may have several vary power functions depending
on the different weight functions are adopted. On the other hand, we may base on a specific weight function to study the
power properties and power dominance's for some various tests. The phenomenon of power dominance on the interesting tests
may vary different weight functions. The crux of using Bayesian concepts to define the level of significance lies in
determining how to decide what is the optimal weight function for the nuisance parameter. Besides, one of the disadvantages
for utilizing Bayesian notions to define the level of significance may require extensive numerical
computations to get the critical values when $p\geq 4$.

\indent As a consequence, for a specific test when the rejected probability under null
hypothesis depends on the nuisance parameter ${\mbSigma}$, its power properties such as
similarity, unbiasedness, and uniformly more powerful phenomenon relative to other tests may
be different if the different ways of obtaining critical values are adopted.
In this paper, we regard the different ways of obtaining critical values for a specific test
as the different kinds of definitions of significance level for that specific test. The way
to give a new definition of the level of significance can be viewed as a way to dig up the
insight information of the null hypothesis. Of course, there are other approaches to define the
level of significance other than using Bayesian notions and maximum principle.

\indent In passing, by virtue of (2.12) and (2.13) we note that the Fisher approach for
$p$-values of the LRT and the UIT are hardly reported since their null distributions depend on
the nuisance parameter. For the problem of testing against the positive orthant space alternative, the
main problem of reporting $p$-value still lies in how to deal with the nuisance parameter.
Several studies for handling the nuisance parameter in reporting $p$-value, which can also be
applied to the present problem, have been well investigated in the literature, we may refer to
Berger et al. [3] for the details. One of the main difficulties is how to find a unified way to
handle the problems. We are afraid that neither the Fisher approach for reporting $p$-value alone
nor the Neyman-Pearson approach alone is well enough to take care well for the problem of testing
against the positive orthant space alternative. Most people think that these two approaches disagree with
each other, however, Lehmann [12] asserted that these two approaches are complementary
rather than contradictory. In passing, we may note that the essential point of Fisher reporting
$p$-value is to find the optimal ways to obtain the insight information under null hypothesis,
and Neyman-Pearson theory essentially shows us how to find the tests that enjoy some optimal
power properties. As such, we may agree with Lehmann that there are essentially only one theory,
not two, for both Fisher and Neyman-Pearson approaches. The spirit lies in the compromise of
Fisher's approach emphasizing on the type I errors and Neyman-Pearson optimal theory without detailed
consideration of power. Therefore we suggest using the unified conservative maximum principle
to define the level of significance $\alpha$ and then adopt the criterion in the decision-theoretic sense
for the power dominance problems.

\vspace{0.3cm}
\def \theequation{3.\arabic{equation}}
\setcounter{equation}{0}
\vspace{0.2cm}

\noindent {\bf 3. A reconciliation of frequentist and Bayesian paradigms}
\vspace{0.2cm}

\indent For testing against the global alternative, Kiefer and Schwartz [11] showed that Hotelling's $T^{2}$-test is a
proper Bayes test. The class of the proper Bayes tests is admissible.

\vspace{0.3cm}
\noindent {\it 3.1 The proper Bayes tests}
\vspace{0.2cm}

\indent The main goal of this section is to see whether the reconciliation of frequentist and Bayesian paradigms can
be achieved for the problem (1.1). First, we establish the necessary conditions for the proper Bayes tests for the
problem (1.1).

\vspace{0.3cm}
\noindent {\bf Proposition 3.1}. {\it Let ${\mathcal A}$ denote the acceptance region of a size-$\alpha$ proper Bayes test
for the problem of testing $H_{0}: {\mbtheta} = {\bf 0}$ against $H_{1}: {\mbtheta} \in {\cal C} \backslash \{{\bf 0}\}$,
where ${\cal C}$ is either the positive orthant space ${\cal O}_{p}^{+}$ or the half-space ${\cal H}_{p}^{*}$. Let
${\mathcal L}$ be any line of support of ${\mathcal A}$, and denote $\Bar{\mathcal A}={\mathcal A}\cup \partial {\mathcal A}$,
$\partial {\mathcal A}$ being the boundary set of ${\mathcal A}$. Then either ${\mathcal L} \subseteq \partial {\mathcal A}$
or ${\mathcal L} \cap \Bar{\mathcal A}=\{{\bf a}\}$, ${\bf a} \in \partial {\mathcal A}$.}

\noindent {\bf Proof}. The density of ${\bf X}_{1}, \ldots, {\bf X}_{n}$ is
\begin{eqnarray}
 \frac{e^{\frac{-1}{2} n{\mbtheta}^{'}{\mbSigma}^{-1}{\mbtheta}}}
   {(2\pi)^{\frac{1}{2}pn}|{\mbSigma}|^{\frac{1}{2}n}}
   \mbox{exp}\left[{n}{\mbtheta}^{'}{\mbSigma}^{-1}{\mXbar}_{n}
    -\frac{1}{2}\mbox{tr}({\mbSigma}^{-1}
    \sum_{j=1}^{n}{\bf X}_{j}{\bf X}^{'}_{j})\right].
\end{eqnarray}
The vector ${\bf y}=({\bf y}^{(1)'},{\bf y}^{(2)'})'$ is composed of ${\bf y}^{(1)}={\mXbar}_{n}$ and
${\bf y}^{(2)}=(v_{11}, 2v_{12},\cdots,2v_{1p}, v_{22},\cdots, v_{pp})$, where
$\left(v_{ij}\right)={\bf V}_{n}=\sum_{j=1}^{n}{\bf X}_{j}{\bf X}^{'}_{j}$.

The vector ${\mbomega}=({\mbomega}^{(1)'},{\mbomega}^{(2)'})'$ is composed of
${\mbomega}^{(1)}=n{\mbSigma}^{-1}{\mbtheta}$ and ${\mbomega}^{(2)}
=-\frac{1}{2}({\sigma}^{11}, {\sigma}^{12},\cdots,{\sigma}^{1p}, {\sigma}^{22}, \cdots,
{\sigma}^{pp})'$, where $({\sigma}^{ij})={\mbSigma}^{-1}$. Write
${\mbSigma}^{-1}={\mbGamma}({\mbomega}^{(2)})$, thus the density of ${\bf Y}$ becomes
\begin{eqnarray}
\frac{|{\mbGamma}({\mbomega}^{(2)})|^{\frac{1}{2}n}\mbox{e}^{-\frac{1}{2n}
{\mbomega}^{(1)'}{\mbGamma}^{-1}({\mbomega}^{(2)}){\mbomega}^{(1)}}}{(2\pi)^{\frac{1}{2}pn}}
  \mbox{exp} \{{\mbomega}^{'}{\bf y}\}
\end{eqnarray}

\indent Since $\mbSigma \in {\cal M}$, thus ${\mbtheta} \in {\cal O}_{p}^{+}$ implies that ${\mbomega}^{(1)} \in {\cal S}$.
Note that every closed convex cone in ${\cal S}$ can be linearly transformed into ${\cal O}_{p}^{+}$ based on the
transformed data with a different covariance matrix. For the sake of simplicity in notation, we may treat every convex cone
in ${\cal S}$ to be ${\cal O}_{p}^{+}$ with transformed data. Similarly, when ${\mbtheta} \in {\cal H}_{p}^{*}$ implies that
${\mbomega}^{(1)}$ represents another half-space with the transformed data. Without loss of generality, we will focus on
the case where ${\mbomega}^{(1)} \in {\cal O}_{p}^{+}$, for the situation of a half-space alternative, which can be handled
similarly. Let $G(\mbomega)$ denote an arbitrary prior distribution function on the set ${\mathcal N}=\{\mbomega|~
{\mbomega}^{(1)}\in {\cal O}_{p}^{+}~\hbox{and}~{\mbGamma}({\mbomega}^{(2)}) \in {\mathcal M}\}$. Then the size-$\alpha$
Bayes test of $H_0$ vs. $H_1$ against this prior distribution rejects $H_0$ if
\begin{eqnarray}
t({\bf y})\triangleq \int_{{\mathcal N}}\mbox{exp}\{{\mbomega}^{(1)'}{\bf y}^{(1)}\}
   dG^{*}({\mbomega}) \geq c,
\end{eqnarray}
where
\begin{eqnarray}
 dG^{*}({\mbomega})=\mbox{e}^{-\frac{1}{2n}{\mbomega}^{(1)'}
  {\mbGamma}^{-1}({\mbomega}^{(2)}){\mbomega}^{(1)}}dG({\mbomega}),
\end{eqnarray}
and c is a constant. Let ${\mathcal L}=\{{\bf z}|~{\bf z}=\rho {\bf u}+(1-\rho){\bf v}, -\infty < \rho < \infty\}$, where
${\bf u}, {\bf v} \in \partial {\mathcal A}$. Next, we study the relationship of the line ${\cal L}$ and the boundary of
acceptance region ${\cal A}$. Let $P^{*}$ be the probability measure induced by $G^{*}$ on the parameter space
${\mathcal N}$.

\indent  Note that (i) if $P^{*}\{{\mbomega}^{'}{\bf (u-v)} \ne 0\}=0$, then
\begin{eqnarray}
t({\bf z})=\int_{{\mathcal N}}\mbox{exp}\{{\mbomega}^{'}{\bf v}+\rho
{\mbomega}^{'}{\bf (u-v)}\}dG^{*}({\mbomega})=c.
\end{eqnarray}
Hence, the line through ${\bf u}$ and ${\bf v}$ is part of the boundary of ${\mathcal A}$.
Namely, ${\mathcal L} \subseteq \partial {\mathcal A}$.

\indent On the other hand, (ii) if $P^{*}\{{\mbomega}^{'}{\bf(u-v)}\ne 0\}\ne 0$, then by the
inequality
\begin{eqnarray}
\mbox{exp}(\rho m_1+(1-\rho) m_2) \leq \rho~ e^{m_1}+(1-\rho)~e^{m_2}, ~0 < \rho~<1,
\end{eqnarray}
for all real $m_1 $ and $m_2$ with equality if and only if $m_1=m_2$, thus we obtain that $t({\bf z}) < c$. Hence
the line ${\mathcal L}$ is not a line of support of ${\mathcal A}$ unless the intersection of ${\mathcal L}$ and
the closure of ${\mathcal A}$ is a single boundary point. Namely, ${\mathcal L} \cap \Bar{\mathcal A}=\{{\bf a}\}$,
${\bf a} \in \partial {\mathcal A}$.
\vspace{0.2cm}

\indent Since the UIT statistics in (2.8) and in (2.10) have a similar structure, hence without loss of generality it
suffices to consider the problem of testing against the positive orthant space alternative.
Let ${\mathcal A}_{U}=\cup_{\emptyset}^{P}{\cal A}_{Ua}$ be the acceptance region of the UIT for testing
$H_{0}: {\mbtheta} = {\bf 0}$ against $H_{1}: {\mbtheta} \in {\cal O}^{+}_{p} \backslash \{{\bf 0}\}$, where
${\cal A}_{Ua}=\{({\mXbar}_{n}, {\bf S}_{n})| ~n{\mXbar}'_{na:a'} {\bf S}^{-1}_{naa:a'}{\mXbar}_{na:a'} \leq u_{\alpha
}\}I_{na},$ {\reg}, with $u_{\alpha}$ be the critical point of the level of significance $\alpha$ for the UIT, which can be
obtained from the right hand sides of (2.16). Consider the case when $p$=2, let ${\cal L}_{a}$ be a line of support of the
set $\Bar{\cal A}_{U} (={\cal A}_{U}\cup \partial {\cal A}_{U})$. Then, (i) for $|a|=2$, we have that ${\cal L}_{a} \cap
\Bar {\cal A}_{Ua}=\{{\bf a}\}, {\bf a} \in \partial {\cal A}_{Ua}$, (ii) for $|a|=1$, we may note that there exists a
line ${\mathcal L}_{a}$ of support of the set ${\cal A}_{Ua}$ such that $({\mathcal L}_{a}\cap I_{na})\subseteq
\partial{\cal A}_{Ua}$ but ${\mathcal L}_{a} \nsubseteq \partial{\cal A}_{Ua}$, since ${\mathcal L}_{a} \cap I_{na}$ is
easily seen to be a half-line, (iii) for $|a|=0$, it is an interior part of the set $\Bar{\cal A}_{U}$. Thus when $p$=2,
the result for the case $|a|=1$ does not satisfy the conditions of Proposition 3.1. The case for testing against the
half-space alternative can be similarly proved, hence we omit the details. Thus, we have the following.

\vspace{0.3cm}
\noindent {\bf Theorem 3.1}. {\it For the problem of testing $H_{0}:{\mbtheta} = {\bf 0}$ against
$H_{1}:{\mbtheta} \in {\cal C}\backslash \{{\bf 0}\}$, where ${\cal C}$ is either the positive orthant space
${\cal O}_{p}^{+}$ or the half-space ${\cal H}_{p}^{*}$, the UIT is not a proper Bayes test.}
\vspace{0.2cm}

\indent Note that for the problem of testing against the positive orthant space alternative or the problem
of testing against the half-space alternative, the LRT is isomorphic to the UIT when the covariance matrix is known
$({\mbSigma}={\mbSigma}_{0}$, with ${\bf S}_{n}$ being replaced by ${\mbSigma}_{0}$) or is partially unknown,
(${\mbSigma}={\sigma}^{2}{\mbSigma}_{0}$, where ${\sigma}^{2}$ is an unknown scalar and ${\mbSigma}_{0}$ is a known
positive definite matrix, with ${\bf S}_{n}$ being replaced by $s^{2}{\mbSigma}_{0}$). Hence, similar arguments as in
the proof of Theorem 3.1, (i) we may conclude that for the problems of testing against restricted alternatives, the LRTs are
not the Bayes tests when the covariance matrix is known or partially unknown. Moreover, for the problems of testing against
restricted alternatives, the LRTs are asymptotic equivalent to the UITs, respectively when the covariance matrix is totally
unknown as the sample size is sufficiently large. Hence, (ii) when the covariance matrix is totally unknown for the problems of
testing against restricted alternatives, the LRTs are not the Bayes tests when sample size is sufficiently large,
(iii) When the covariance matrix is totally unknown and the sample size is fixed, let ${\mathcal A}_{L}
=\cup_{\emptyset}^{P}{\cal A}_{La}$ be the acceptance region of the LRT for the problem of testing against the positive
orthant space alternative, $H_{0}: {\mbtheta} = {\bf 0}$ against $H_{1}: {\mbtheta} \in {\cal O}^{+}_{p} \backslash
\{{\bf 0}\}$, where ${\cal A}_{La}=\{({\mXbar}_{n}, {\bf S}_{n})|~{\frac{n{\mXbar}'_{na:a'}{\bf S}_{naa:a'}^{-1}
{\mXbar}_{na:a'}}{1+n{\mXbar}'_{na'} {\bf S}_{na'a'}^ {-1} {\mXbar}_{na'}}} \leq l_{\alpha} \}I_{na}$ with $l_{\alpha}$
being the critical point of the level of significance $\alpha$ for the LRT. The $l_{\alpha}$ can be obtained from the
right hand sides of (2.15). Consider the case when $p$=2, (i) for $|a|=2$, there exists a line ${\cal L}_{a}$ of support of
${\cal A}_{L}$ so that there are more than one intersection point $({\cal L}_{a}\cap \Bar{\cal A}_{L})$, (ii) for $|a|=1$,
we may have that ${\cal L}_{a}\cap \Bar{\cal A}_{L\alpha}=\{{\bf a}\}, {\bf a} \in \partial {\cal A}_{L\alpha}$, (iii) for
$|a|=0$, it is an interior part of the set $\Bar{\cal A}_{L\alpha}$. Thus when $p$=2, the result for the case $|a|=2$
does not satisfy the conditions of Proposition 3.1. The case for testing against the half-space alternative can be
similarly proved, hence we omit the details. Thus, we have the following.

\vspace{0.3cm}
\noindent {\bf Theorem 3.2}. {\it For the problem of testing $H_{0}:{\mbtheta} = {\bf 0}$ against
$H_{1}:{\mbtheta} \in {\cal C}\backslash \{{\bf 0}\}$, where ${\cal C}$ is either the positive orthant space
${\cal O}_{p}^{+}$ or the half-space ${\cal H}_{p}^{*}$, the LRT is not a proper Bayes test.}
\vspace{0.2cm}

\indent Next, consider for the finite union-intersection test (FUIT) of Roy et al. [22] based on the one-sided
coordinated wise Student t-tests. Define ${\bf S}_{n}=(s_{nij})$ as in (2.1) and let
\begin{eqnarray}
t_j=\frac{{\sqrt {n}}\bar X_{nj}}{\sqrt {s_{njj}}},~ j=1, \cdots, p.
\end{eqnarray}
Corresponding to a given significance level $\alpha$, define
\begin{eqnarray}
{\alpha}^{*}:~p{\alpha}^{*}=\alpha.
\end{eqnarray}
Let then $t_{n-1,{\alpha}^{*}}$ be the upper
$100{\alpha}^{*}$\% point of the Student t-distribution with $n-1$ degrees of freedom.
Consider the critical region ${\cal W}_j=\{t_j|~t_j\geq t_{n-1,{\alpha}^{*}}\}$ for
$~j=1,\cdots, p$. Then the critical region of the FUIT is
\begin{eqnarray}
{\cal W}^{*}=\cup_{j=1}^{p}{\cal }{\cal W}_j
\end{eqnarray}
and the acceptance region is
\begin{eqnarray}
{\cal A}^{*}=\cap_{j=1}^{p}\Bar{\cal A}_j,
\end{eqnarray}
where ${\cal A}_j=R\setminus{\cal W}_j$ and $\Bar{\cal A}_j$ denote the closure of ${\cal A}_j, j=1, \cdots, p.$

\indent Let ${\cal A}=\cap_{j=1}^{p}{\cal A}_j$ and ${\mathcal L}_{a}$ be a line of support of ${\mathcal A}$. Consider the
case when $p$=2, (i) for $|a|=2$ and $|a|=1$, then we have that $({\cal L}_{a}\cap I_{na}) \subseteq {\cal A}^{*}$, but
${\cal L}_{a} \nsubseteq {\cal A}$, (ii) for $|a|=0$, it is an interior part of the set ${\cal A}^{*}$. Thus when $p$=2,
the results for the cases $|a|=2$ and $|a|=1$ do not satisfy the conditions of Proposition 3.1. The case for testing
against the half-space alternative can be similarly proved, hence we omit the details.Thus, we have the following.

\vspace{0.3cm}
\noindent {\bf Theorem 3.3}. {\it For the problem of testing $H_{0}:{\mbtheta} = {\bf 0}$ against
$H_{1}:{\mbtheta} \in {\cal C}\backslash \{{\bf 0}\}$, where ${\cal C}$ is either the positive orthant space
${\cal O}_{p}^{+}$ or the half-space ${\cal H}_{p}^{*}$, the FUIT is not a proper Bayes test.}
\vspace{0.2cm}

\indent For testing against the global alternative, the LRT and the UIT are equivalent to Hotelling's $T^{2}$-test, which is a proper
Bayes test. However, for the problem of testing against restricted alternatives, Theorems 3.1-3.3 tell us that the UIT, the
LRT, and the FUIT are not the proper Bayes tests. The class of admissible tests is essentially larger than the class of the
class of proper Bayes tests, as it will be demonstrated later on. Obviously, for the problem (1.1), Proposition 3.1
provides us with an easy way to examine whether the test is a proper Bayes test. We may note that the explicit forms of
the proper Bayes tests are hardly obtained due to the difficulty of integration over restricted parameter spaces. Hence,
for the problems of testing against restricted alternatives we will study the other non-proper Bayesian tests regarding
the non-informative priors for the nuisance parameter ${\mbSigma}$.

\vspace{0.3cm}
\noindent {\it 3.2 The integrated LRT and the integrated UIT}
\vspace{0.2cm}

\indent Efron [10] gave the 1996 Fisher lecture entitled "R. A. Fisher in the 21 Century". Professor Rob Kass, one of the
followed discussants, commented that ``the situation in statistical inference is beautifully peaceful and compelling for
one-parameter problems. When we go to the multiparameter world, however, the hope dims out not only for a reconciliation
of Bayesian and frequentist paradigms, but for any satisfactory, unified approach in either a frequentist or Bayesian
framework." According to his comments and the Efron's lecture, we may realize that reconciliation of Bayesian and
frequentist paradugms is an important and fundamental issue in theoretical statistics.

\indent The likelihood integrated (Berger et al. [3]) with respect to a relevant Haar measure for the nuisance covariance
matrix while remaining neutral with respect to the parameter of interest, which fiducial inference was intended to be,
can be viewed in a objective Bayesian light.

\indent Riemannian geometry yields an invariant volume element $dv$ on ${\mathcal M}$ which is of the form
$dv=|{\mbSigma}|^{-(p+1)/2}(d{\mbSigma})$, where $(d{\mbSigma})=\prod_{1 \leq j \leq i  \leq p} d{\sigma}_{ij}$
with $d{\sigma}_{ij}$ being the Lebesgue measure on $R$. We may naturally adopt this relevant Haar measure $dv$ as the
underline measure of the a prior of $\mbSigma$. The likelihood integrated with respect to this relevant Haar measure
$dv$ for the nuisance covariance matrix can be viewed as in a Bayesian light. It is not hard to note that for testing
against global alternative, Hotelling's $T^{2}$-test is the version of integrated LRT (and UIT), see Berger et al. [3].

\indent Write
\begin{eqnarray}
   g({\mbSigma}) \propto  ~|{\mbSigma}|^{-\frac {1}{ 2}(p+1)} ~d\nu({\mbSigma}),
            ~ {\mbSigma} \in {\mathcal M}.
\end{eqnarray}
First, consider the unique Haar measure in ${\mbSigma}$ over the set ${\cal M}$, i.e., $d\nu({\mbSigma})=d{\mbSigma}$.
Then the marginal density function of ${\mXbar}_{n}$ and ${\bf S}_{n}$ is given by
\begin{eqnarray}
  P_{1}({\mXbar}_{n}, {\bf S}_{n}|~{\mbtheta}) \propto ~|{\bf S}_{n}|^{\frac{1}{2}(n-p-2)}
  |{\bf S}_{n}+n({\mXbar}_{n}-{\mbtheta})({\mXbar}_{n}-{\mbtheta})'|^{-\frac{1}{2}(n-1)},
\end{eqnarray}
namely,
\begin{eqnarray}
  P_{1}({\mXbar}_{n}, {\bf S}_{n}|~{\mbtheta}) \propto ~|{\bf S}_{n}|^{-\frac{1}{2}(p+1)}
   [1+n({\mXbar}_{n}-{\mbtheta})'{\bf S}_{n}^{-1}
   ({\mXbar}_{n}-{\mbtheta})]^{-\frac{1}{2}(n-1)}.
\end{eqnarray}

\indent Based on the density function in (3.13), we then consider the integrated likelihood ratio
\begin{eqnarray}
&~&\frac{\sup_{\mbtheta \in {\cal O}_{p}^{+}}P_{1}({\mXbar}_{n}, {\bf S}_{n}|~{\mbtheta})}
{P_{1}({{\mXbar}_{n}, {\bf S}_{n}|~{\bf 0})}} \\
  & = & \left[ \frac{1+n{\mXbar}_{n}'{\bf S}_{n}^{-1}{\mXbar}_{n}}
    {1+\inf_{\mbtheta \in {\cal O}^{+}_{p}}n({\mXbar}_{n}-{\mbtheta})'{\bf S}_{n}^{-1}
     ({\mXbar}_{n}-{\mbtheta})} \right ]^{\frac{1}{2}(n-1)} \nonumber \\
 & = & \left [1+|\!| {\pi}_{{\bf S}_{n}}(n^{1/2} {\mXbar}_{n}; {\cal O}^{+}_{p})
    |\!|_{{\bf S}_{n}}^{2} \{ 1 + |\!| n^{1/2} {\mXbar}_{n} - {\pi}_{{\bf S}_{n}}
     (n^{1/2} {\mXbar}_{n}; {\cal O}^{+}_{p} ) |\!|_{{\bf S}_{n}}^{2} \}^{-1}\right
     ]^{\frac{1}{2}(n-1)} \nonumber \\
  & = & (1+ L_{n})^{\frac{1}{2}(n-1)}, \nonumber
\end{eqnarray}
which is the monotone function of LRT. Thus, we may obtain that the LRT is the integrated LRT for the problem of testing
against the positive orthant space alternative.

\indent  As for the integrated UIT, note that for any ${\bf b} \in R^{p}\backslash\{{\bf 0}\}$, the following
density function can be obtained from (3.13)
\begin{eqnarray}
  P_{2}({\bf b}^{'}{\mXbar}_{n}, {\bf b}^{'}{\bf S}_{n}{\bf b}|~{\mbtheta}) &\propto& ~|{\bf b}'{\bf S}_{n}
  {\bf b}|^{-\frac{1}{2}(p+1)} [1+n {\bf b}'({\mXbar}_{n}-{\mbtheta})
  ({\mXbar}_{n}-{\mbtheta})'{\bf b}/{\bf b}'{\bf S}_{n}{\bf b}]^{-\frac{1}{2}(n-1)}.
\end{eqnarray}
Based on the density function in (3.15), for each ${\bf b} \in R^{p}\backslash\{{\bf 0}\}$, we have
\begin{eqnarray}
&~&\frac{\sup_{\mbtheta \in {\cal O}_{p}^{+}}P_{2}({\bf b}^{'}{\mXbar}_{n}, {\bf b}^{'}{\bf S}_{n}{\bf b}|~{\mbtheta})}
{P_{2}({\bf b}^{'}{\mXbar}_{n}, {\bf b}^{'}{\bf S}_{n}{\bf b}|~{\bf 0})} \\
  & = & \left[ \frac{1+n {\bf b}'{\mXbar}_{n}{\mXbar}'_{n}{\bf b}/
    {\bf b}'{\bf S}_{n} {\bf b}}{1 + \inf_{\mbtheta \in {\cal O}_{p}^{+}}
    n({\mXbar}_{n}-{\mbtheta})'{\bf b} {\bf b}'({\mXbar}_{n}-{\mbtheta})/
    {\bf b}' {\bf S}_{n}{\bf b}} \right ]^{\frac{1}{2}(n-1)}. \nonumber
\end{eqnarray}
Partition ${\bf b}'=({\bf b}'_{a}, {\bf b}'_{a'})$ the same as in (2.3) and for each
$a$, \reg, let ${\bf G}_{a}= \left[
\begin{array}{cc}
    {\bf I}_{a}& {\bf -S}_{naa'}{\bf S}^{-1}_{na'a'}\\
    {\bf 0}& {\bf I}_{a'}
\end{array}
  \right].$
Further multiply ${\bf G}_{a}$ on the left of ${\mbtheta}$ and ${\mXbar}_{n}$
and its inverse on the right of ${\bf b}^{'}$. Next by using Lemma 1.1 of N\"{u}esch
[17] and choosing ${\bf b}_{a'}= - {\bf S}_{na'a'}{\bf S}_{na'a}{\bf b}_{a}$,
then after some simplifications the $ r.h.s.$ of (3.16) reduces to
\begin{eqnarray}
\left[ 1+ \sum_{\mbox{\reg}} \left\{ {\frac{n({\bf b}'_{a}
{\mXbar}_{na:a'})^{2}}
{{\bf b}'_{a}{\bf S}_{naa:a'}{\bf b}_{a}}} \right\}I_{na} \right]^{\frac{1}{2}(n-1)}.
\end{eqnarray}
Since (3.17) holds for all ${\bf b}_{a} \neq {\bf 0}$, \reg, further note that
\begin{eqnarray}
&~&\sum_{\mbox{\reg}}\sup_{{\bf b}_{a} \neq {\bf 0}}\left\{ {\frac{n({\bf b}'_{a}
        {\mXbar}_{na:a'})^{2}}{{\bf b}'_{a}{\bf S}_{naa:a'}{\bf b}_{a}}} \right\}I_{na}\\
& = &  \sum_{\mbox{\reg}}\{n{\mXbar}'_{na:a'}{\bf S}_{naa:a'}^
         {-1}{\mXbar}_{na:a'}\}I_{na} \nonumber \\
& = & U_{n}, \nonumber
\end{eqnarray}
thus the UIT is the integrated UIT for the problem of testing against the positive orthant space alternative.

\indent Next, we consider another improper prior of $\mbSigma$, which is called the Bernardo reference prior (Chang and
Eaves [8]), of the following form
\begin{eqnarray}
  g_{1}({\mbSigma}) \propto ~|{\mbSigma}|^{-\frac{1}{2}(p+1)}
     |{\bf I}+{\mbSigma}*{\mbSigma}^{-1}|^{-\frac{1}{2}}~d{\mbSigma},
     ~{\mbSigma}\in {\mathcal M},
\end{eqnarray}
where ${\bf C}*{\bf D}$ denotes the Hadamard product of the matrices ${\bf C}$ and ${\bf D}$. Namely, $d\nu({\mbSigma})
=|{\bf I}+{\mbSigma}*{\mbSigma}^{-1}|^{-\frac{1}{2}}~d{\mbSigma}$ in (3.19). Then after some manipulation, the marginal
density function is
\begin{eqnarray}
  h({\bf W}_{n})= c \int_{{\mbSigma}\in {\mathcal M}}
               |{\mbSigma}|^{-\frac{1}{2}(p+1)}|{\bf I}+
         {\mbSigma}*{\mbSigma}^{-1}|^{-\frac{1}{2}}~
          \mbox{e}^{{-\frac{1}{2}}\mbox{tr}({\mbSigma}^{-1}
          {\bf W}_{n})}~d{\mbSigma},
\end{eqnarray}
where ${\bf W}_{n}={\bf S}_{n}+n({\mXbar}_{n}-{\mbtheta})({\mXbar}_{n}-{\mbtheta})'$ and $c$ is the normalizing constant.
By the facts that $d|{\bf W}_{n}|= |{\bf W}_{n}|\mbox{tr}({\bf W}^{-1}_{n})(d{\bf W}_{n})$ and
$d\mbox{tr}({\mbSigma}^{-1}{\bf W}_{n})=\mbox{tr}({\mbSigma}^{-1})(d{\bf W}_{n})$, where $(d{\bf W}_{n})$ denotes the
exterior differential form of ${\bf W}_{n}$. Then
\begin{eqnarray}
\frac{dh({\bf W}_{n})}{d|{\bf W}_{n}|}
    &=& {\frac{-c}{2}}\int_{{\mbSigma}\in {\mathcal M}}|{\bf W}_{n}|^{-1}
         (\mbox{tr}{\bf W}^{-1}_{n})^{-1}\mbox{tr}({\mbSigma}^{-1})
         |{\mbSigma}|^{-\frac{1}{2}(p+1)}  \\
     & & \hspace{3cm} \times |{\bf I}+{\mbSigma}*{\mbSigma}^{-1}|^{-\frac{1}{2}}
         \mbox{e}^{{-\frac{1}{2}}\mbox{tr}({\mbSigma}^{-1}{\bf W}_{n})}
         ~d{\mbSigma} \nonumber  \\
     &< &  0 , \nonumber
\end{eqnarray}
hence $h({\bf W}_{n})$ is strictly decreasing in $|{\bf W}_{n}|$. Thus the corresponding integrated LRT and integrated UIT
for the Bernardo reference prior can be obtained as those of the Haar measure case, respectively.

\indent Finally, we use a proper prior of ${\mbSigma}$ to check how is the difference to those improper priors mentioned
above. Consider the case if ${\mbSigma}$ is assigned a proper prior distribution $W^{-1}_{p}(\mbGamma, m)$, where
$W^{-1}_{p}(\mbGamma, m)$ denoting the $p$-dimensional inverted Wishart distribution with $m$ degrees of freedom and
expectation $m{\mbGamma}$ which is positive definite. Then by Theorem 7.7.2 of Anderson [1], the marginal density is
\begin{eqnarray}
  h_{1}({\mXbar}_{n}, {\bf S}_{n}|~{\mbtheta}) & \propto & ~|{\mbGamma}|^{\frac{m}{2}}
      |{\bf S}_{n}|^{\frac{1}{2}(n-p-2)}|{\bf S}_{n}+n({\mXbar}_{n}
      -{\mbtheta})({\mXbar}_{n} -{\mbtheta})^{'}+{\mbGamma}|^{-\frac{1}{2}(n+m-1)}   \\
    & = &  |{\mbGamma}|^{\frac{m}{2}}|{\bf S}_{n}|^{\frac{1}{2}(n-p-2)}
         |{\bf W}_{n}+{\mbGamma}|^{-\frac{1}{2}(n+m-1)}, \nonumber
\end{eqnarray}
which is a monotone function of $|{\bf W}_{n}|$. Thus, we may conclude that the corresponding integrated LRTs and
integrated UITs with respect to the Haar measure and the Bernard reference prior are the same as that of ${\mbSigma}$
being a proper prior distribution $W^{-1}_{p}(\mbGamma, m)$.

\indent Similarly, for the problem of testing against the half-space alternative ${\cal H}_{p}^{*}$ by (3.14) and (3.17) the
$L_{n}^{*}$ and $U_{n}^{*}$ can be exactly obtained with respect to those non-informative prior differentials and the
inverted Wishart distribution mentioned above over the space ${\mbSigma}\in {\mathcal M}$.

\indent Moreover, for the problem (1.1) the corresponding integrated LRTs and integrated UITs also have the same forms as
those of the corresponding LRTs and UITs, respectively. Based on the discussions, we may conclude in the following.

\vspace{0.3cm}
\noindent {\bf Theorem 3.4}. {\it For the problem of testing $H_{0}:{\mbtheta} = {\bf 0}$ against $H_{1}:{\mbtheta}
\in {\cal C}\backslash \{{\bf 0}\}$, where ${\cal C}$ is either the positive orthant space ${\cal O}_{p}^{+}$ or the
half-space ${\cal H}_{p}^{*}=\{{\mbtheta}|~{\theta}_{p} \geq 0 \}$, then both the LRT and the UIT are the versions of
the integrated LRT and the integrated UIT, respectively.}
\vspace{0.2cm}

\indent Hence, we may note that Theorem 3.4 results in the reconciliation of frequentist and Bayesian paradigms via the
fiducial inference for the problem (1.1).  For problem (1.1), we continue to study whether any satisfactory, unified
approach or framework for the LRT or the UIT exists.

\vspace{0.3cm}
\def \theequation{4.\arabic{equation}}
\setcounter{equation}{0}

\vspace{0.3cm}
\noindent {\bf 4. Wald's decision theory}
\vspace{0.2cm}

\indent For the problem of testing against the positive orthant space alternative, the difficulty is that both the
distribution functions of the LRT and the UIT under the null hypothesis are dependent on the nuisance parameter. In the
literature, many authors heavily rely on the Neyman-Pearson optimal theory to construct new tests so that the new tests
dominate the corresponding LRT in the sense of being uniformly more powerful. For the details, we refer to Berger [4],
Menendezand Salvador [15], and Menendez et al. [16] and the references therein. Sen and Tsai [24] adopted Stein's concept
[26] to establish the two-stage LRT, which enjoys some optimal properties, however, the cost is the complicated distribution
function of the test statistic. The method proposed by Sen and Tsai [24] is fundamentally based on the Neyman-Pearson
optimal theory.

\indent For testing against a global shift alternative, Hotelling's $T^{2}$-test is uniformly the most powerful invariant
(UMPI), and hence, is also admissible (Simika [25]). Stein [27] established the admissibility of Hotelling's
$T^{2}$-test by using the exponential structure of the parameter space. The UMPI character, as well as the admissibility
of Hotelling's $T^{2}$-test may not generally hold for restricted alternatives, such as $H_{1{\cal O}^{+}}$ in (1.1).
The affine-invariance structure of the parameter space $\mathrm {\Theta} =\{{\mbtheta}\in R^{p}\}$ does not hold for
$H_{1{\cal O}^{+}}$, and hence, when ${\mbSigma}$ is arbitrary p.d., restriction to invariant tests makes little sense.
As such, it is conjectured, though not formally established, that possibly some other non-(affine) invariant tests
dominating Hotelling's $T^{2}$-tests, and hence, the latter is inadmissible. Our objective is to establish that
for testing $H_0$ vs $H_{1{\cal O}^{+}}$, Hotelling's $T^{2}$-test is inadmissible. For testing $H_0$ against
$H_{1{\cal O}^{+}}$, the set of proper Bayes tests and their weak limits only constitute a proper subset of an essentially
complete class of tests (Marden [14]). Marden's minimal complete classes consist of proper Bayes tests and the tests with
convex and decreasing acceptance regions for density functions more general than exponential family. For the problem
(1.1), Proposition 3.1 tells us that the UITs are not the proper Bayes. Furthermore, the acceptance regions of the UITs may
not be decreasing, though they can be demonstrated to be the convex sets. Therefore, the conditions of Marden's minimal
complete classes of tests are also insufficient for the UITs. For problem (1.1), it looks hopeless that the UITs are
admissible.

\vspace{0.3cm}
\noindent {\it 4.1 Eaton's essentially complete class of tests}
\vspace{0.2cm}

\indent We evaluate Eaton's [9] fundamental result regarding an essentially complete class of test functions for testing
against restricted alternatives, when the underlying density belongs to an exponential family, as in our current context.
Let $\Phi$ be Eaton's essentially complete class of tests. This means for any test ${\varphi}^{*}\notin {\Phi}$ there
exists a test ${\varphi} \in {\Phi}$ such that ${\varphi}$ is at least as good as ${\varphi}^{*}$.

\vspace{0.3cm}
\noindent {\bf Theorem 4.1.} {\it For the problem of testing $H_{0}: {\mbtheta} = {\bf 0}$ against $H_{1{\cal O}^{+}}:
{\mbtheta} \in {\cal O}_{p}^{+} \backslash \{{\bf 0}\}$, Hotelling's $T^{2}$-test is inadmissible.}

\noindent {\bf Proof}. The present testing hypothesis problem is invariant under the group of transformations of positive
diagonal matrices, hence, for simplicity, $\mbSigma$ may be treated as the correlation matrix. Following Eaton [9],
we define
\begin{eqnarray}
 \Omega_{1}=\{{\mbSigma}^{-1}{\mbtheta}|~{\mbtheta} \in {\cal O}_{p}^{+}\}
 \backslash \{{\bf 0}\}\in {\cal S}\backslash \{{\bf 0}\}.
\end{eqnarray}
Since any closed convex cone can be a suitably linearly transformed to ${\cal O}_P^{+}$ with a different covariance
matrix, and the transformed data ({\bf Y}). It is important to observe that the test statistics based on the transformed
data have the same format as those based on the original data. The statistics are based on the original data in the
probability space $({\mXbar}_{n},{\bf S}_{n})$ and a new transformed one under the new probability space $({\mYbar}_{n},
{\bf S}_{ny})$ are exactly the same with probability one. For the sake of simplicity in notations, we will continue using
the original data notation, namely directly treating the cone as ${\cal O}_{p}^{+}$, i.e., $\Omega_{1}={\cal O}^{+}_{p}
\backslash \{{\bf 0}\}$.

\indent Let ${\cal V}\subseteq R^{p}$ be the smallest closed convex cone containing $\Omega_{1}$. Then the dual cone of
${\cal V}$ is defined as
\begin{eqnarray}
 {\cal V}^{-}=\{{\bf w}|~<{\bf w}, {\bf x}>_{\mbSigma} \leq 0,~\forall~ {\bf x}\in {\cal V}\}.
\end{eqnarray}
Then, we have that ${\cal V}={\cal O}_{p}^{+}$, and its dual cone is
\begin{eqnarray}
{\cal V}^{-}&=&\{{\bf w}|~{\bf x}^{'}{\mbSigma}^{-1}{\bf w} \leq 0, ~\forall~ {\bf x}\in {\cal V}\}. \\ \nonumber
            &=&\{{\bf w}|~{\mbSigma}^{-1}{\bf w} \leq {\bf 0} \},  \nonumber
\end{eqnarray}
which is an unbounded set. Note that a test that is a member of Eaton's [9] essentially complete class of tests
can be established by showing the acceptance region contains the subspace $\mathcal {V}^{-}$.

\indent The acceptance region of Hotelling's $T^{2}$-test is given by
\begin{eqnarray}
\mathcal {A}_{T^{2},t^{2}_{\alpha}}=\{({\mXbar}_{n},{\bf S}_{n})|~T^{2} \leq {t^{2}_{\alpha}}\},
\end{eqnarray}
where ${t^{2}_{\alpha}}$ is the upper $100\alpha$\% point of the null hypothesis distribution of $T^{2}$ (which is linked
to a F-distribution). Since $\mathcal {A}_{T^{2},t^{2}_{\alpha}}$ is a hyper-ellipsoidal set with origin ${\bf 0}$, it is
bounded, while $\mathcal {V}^{-}$, as shown above, is an unbounded set. Therefore, Eaton's [9] condition is not tenable,
Hotelling's $T^{2}$-test is not a member of Eaton's essentially complete class of tests. Thus, the theorem follows.

\vspace{0.3cm}
\indent Similarly, for the problem of testing $H_{0}$ vs $H_{1{\cal H}^{*}}$. When ${\Omega}_{1}={\cal H}_{p}^{*} \backslash
\{{\bf 0}\}$, then we have that ${\cal V}={\cal H}_{p}^{*}$. Hence its dual cone ${\cal V}^{-}=\{{\mbtheta}\in
R^{p}|~{\theta}_{p}=0\}$ being the hyperplane, is an unbounded set. Thus, Hotelling's $T^{2}$ for testing against the
half-space alternative is not a member of Eaton's essentially complete class. Thus we have the following.

\vspace{0.3cm}
\noindent {\bf Theorem 4.2}. {\it For the problem of testing $H_{0}: {\mbtheta}={\bf 0}$ against $H_{1{\cal H}^{*}}:
{\mbtheta} \in {\cal H}_{p}^{*} \backslash \{{\bf 0}\}$, Hotelling's $T^{2}$-test is inadmissible.}
\vspace{0.2cm}

\indent For the problem of testing $H_{0}: {\mbtheta} = {\bf 0}$ against $H_{1{\cal O}^{+}}: {\mbtheta} \in {\cal O}_{p}^{+}
\backslash \{{\bf 0}\}$, the acceptance region of UIT is
\begin{eqnarray}
  {\cal A}_{U}=\left\{({\mXbar}_{n}, {\bf S}_{n})|~\sum_{\mbox{\reg}}\{n{\mXbar}'_{na:a'}
  {\bf S}_{naa:a'}^{-1}{\mXbar}_{na:a'}\}I_{na}\leq k \right \},
\end{eqnarray}
for a suitable $k$. Note that
\begin{eqnarray}
  {\cal A}_{{U}}=\cup_{\emptyset}^{P}{\cal A}_{Ua},
\end{eqnarray}
where
\begin{eqnarray}
  {\cal A}_{Ua}=\{({\mXbar}_{n}, {\bf S}_{n})|~n{\mXbar}'_{na:a'}
    {\bf S}_{naa:a'}^{-1}{\mXbar}_{na:a'}\leq k,
   ~{\mXbar}_{na:a'}>{\bf 0}, {\bf S}_{na'a'}^{-1}{\mXbar}_{na'}\leq {\bf 0}\},
\end{eqnarray}
$\emptyset \subseteq a \subseteq P$.

\indent For testing $H_0$ against$H_{1{\cal O}^{+}}$, though the UIT (Sen and Tsai [24]) is not the proper Bayes we will
show that it is the member of Eaton's [9] essentially complete class. Corresponding to a preassigned $\alpha~(0
< \alpha < 1)$, let $c_{\alpha}$ be the critical level, obtained by equating the right hand side of (2.16) to $\alpha$,
thus the UIT is a size-$\alpha$ test for $H_0$ vs $H_{1{\cal O}^{+}}$. Let ${\cal A}_{U}$ be the acceptance region formed
by letting $U_{n} \leq c_{\alpha}$ in (2.8). Let ${\mathcal V}^{-}_{n}={\cal A}_{U\emptyset}=\{({\mXbar}_{n}, {\bf S}_{n})|~
{\bf S}_{n}^{-1}{\mXbar}_{n} \leq {\bf 0}\}$, since ${\bf S}_{n}$ is positive definite with probability one, we then have
${\mathcal V}^{-}_n \to {\mathcal V}^{-}$. Thus, ${\mathcal V}^{-}_n \to  {\mathcal V}^{-}  \subseteq {\mathcal A}_{U}
-{\bf a}^{0}$ for each ${\bf a}^{0}\in \partial {\mathcal A}$.

\indent Similarly, for the FUIT (Tsai and Sen [30]) we may note that its acceptance set ${\cal A}^{*}$ [ see (3.10)] is a
closed convex set. Note that ${\bf A}=(a_{ij})$ such that $a_{ij} \leq 0$ for all $i \ne j$ is an $M$-matrix if and only
if ${\bf A}^{-1}=(a^{ij}) ~\geq 0~ \forall ~i \ne j$. When ${\mbSigma}$ is an $M$-matrix, then ${\mathcal V}^{-}
=\{{\bf x}\in R^{p}|~{\bf x} \leq {\bf 0}\} \subseteq {\cal A}^{*}$ as long as $t_{n-1, 1-{\alpha}^{*}}$ $\geq 0$,
${\alpha}^{*} \leq \frac{1}{p}$ (or $\alpha \leq 1$). Thus, when $\mbSigma$ is an $M$-matrix the FUIT is a size-$\alpha$
test for $H_0$ vs $H_{1{\cal O}^{+}}$ and ${\mathcal {V}^{-}} \subseteq {\cal A}^{*}$. Thus, when $\mbSigma$ is an
$M$-matrix, the FUIT is also the members of Eaton's essentially complete class for $H_0$ against$H_{1{\cal O}^{+}}$.

\vspace{0.3cm}
\noindent {\bf Theorem 4.3}.  {\it For the problem of testing $H_{0}: {\mbtheta} = {\bf 0}$ against $H_{1{\cal O}^{+}}:
{\mbtheta} \in {\cal O}_{p}^{+} \backslash \{{\bf 0}\}$, the UIT belongs to Eaton's essentially complete class of tests,
so does the FUIT when ${\mbSigma}$ is an $M$-matrix.}
\vspace{0.2cm}

\indent For the problem of testing $H_{0}: {\mbtheta} = {\bf 0}$ against
$H_{1 {\cal H}^{*}}: {\mbtheta} \in {\cal H}^{*}_{p} \backslash \{{\bf 0}\}$, the acceptance region of UIT is
\begin{eqnarray}
  {\cal A}_{{U}^{*}}=\left\{({\mXbar}_{n}, {\bf S}_{n})|~ \sum_{\mbox{\rag}}
   \{ n{\mXbar}'_{na:a'}{\bf S}^{-1}_{naa:a'}{\mXbar}_{na:a'}\}
   I^{*}_{na}\leq k_{1} \right \},
\end{eqnarray}
where $P_{1}=P-\{p\}$ for a suitable $k_{1}$. Note that
\begin{eqnarray}
  {\cal A}_{{U}^{*}}={\cal A}_{U^{*}P} \cup {\cal A}_{U^{*}P_{1}},
\end{eqnarray}
where
\begin{eqnarray}
  {\cal A}_{U^{*}P}=\{({\mXbar}_{n}, {\bf S}_{n})|~n{\mXbar}'_{n}{\bf S}^{-1}_{n}
  {\mXbar}_{n}\leq k_{1},~{\Bar X}_{np} \geq  0\},
\end{eqnarray}
and
\begin{eqnarray}
  {\cal A}_{U^{*}P_{1}}=\{({\mXbar}_{n}, {\bf S}_{n})|~n{\mXbar}'_{nP_{1}: p}
   {\bf S}^{-1}_{nP_{1}P_{1}:p}{\mXbar}_{nP_{1}:p}
    \leq k_{1},~ \Bar{X}_{np} \leq 0 \}
\end{eqnarray}
respectively.

\indent For testing $H_{0}$ vs $H_{1{\cal H}^{*}}$, let $c^{*}_{\alpha}$ be the critical level, obtained by equating the
right hand side of (2.20) to $\alpha$, thus the UIT is a size-$\alpha$ test for $H_0$ vs $H_{1{\cal H}^{*}}$. Let
${\cal A}_{U^{*}}$ be the acceptance region formed by letting $U^{*}_{n} \leq c^{*}_{\alpha}$ in (2.10). it is easy to
note that ${\cal V}^{-}=\{{\mbtheta} \in R^{p}~|~\theta_{p}=0\} \subseteq {\cal A}_{{U}^{*}}$.

\indent Similarly, for the FUIT we may note that its acceptance set ${\cal A}^{h}$ which is the intersection of the
acceptance regions of (p-1) two-sided tests and a one-sided acceptance region, i.e.,${\cal A}^{h}=\{t_{j}|~t_{n-1,
\alpha^{*}/2} \leq t_{j} \leq t_{n-1, 1-\alpha^{*}/2}, j=1, 2, \ldots, p-1$, and $t_{p} \leq t_{n-1, 1-\alpha^{*}}\}$.
It is a cylinder-type convex cone and ${\mathcal V}^{-} \subseteq {\cal A}^{h}$ as long as $t_{n-1, {\alpha}^{*}/2}$
$\leq 0$ and $t_{n-1, 1-{\alpha}^{*}/2} \geq 0$, ${\alpha}^{*} \leq \frac{1}{p}$ (or $\alpha \leq 1$). Thus, the FUIT is
a size-$\alpha$ test for $H_0$ vs $H_{1{\cal H}^{*}}$ and ${\mathcal {V}^{-}} \subseteq {\cal A}^{h}$.

\vspace{0.3cm}
\noindent {\bf Theorem 4.4}.  {\it For the problem of testing $H_{0}: {\mbtheta} = {\bf 0}$ against $H_{1{\cal H}^{*}}:
{\mbtheta} \in {\cal H}_{p}^{*} \backslash \{{\bf 0}\}$, both the UIT and the FUIT belong to Eaton's essentially
complete class of tests.}
\vspace{0.2cm}

\indent By virtues of Theorems 4.1 and 4.2 it is interesting to see whether Hotelling's $T^{2}$-test is dominated by
the UIT. For the case that both Hotelling's $T^{2}$-test and the UIT are unbiased, Tsai and Sen [30] provided an
affirmed answer for it when $\mbSigma={\sigma}^{2} {\mbDelta}$, where ${\sigma}^{2}$ is an unknown scalar parameter
and ${\mbDelta}$ is a known $M$-matrix. However, it is not true when ${\mbSigma}$ is totally unknown due to the fact that
Hotelling's $T^{2}$-test is unbiased while the UIT is not for the problem of testing against the positive orthant space
alternative. Theorems 4.3 and 4.4 do not guarantee that the corresponding UITs are admissible for the problem (1.1).

\vspace{0.3cm}
\noindent {\it 4.2 The Birnbaum-Stein method}
\vspace{0.2cm}

\indent It is interesting to determine whether the UITs are admissible in the decision-theoretic sense (Wald [31]) for
the problem (1.1). Our problems fall under the framework of the exponential family, Theorem 5.6.5 of Anderson [1], which
states that for the conditions for tests to be $d$-admissible for the context of testing against restricted alternatives,
which can be viewed as a generalization of the work by Birnbaum [6] and Stein [27].

\indent For the definitions of $\alpha$-admissible and $d$-admissible, we refer to the book  of Lehmann ([13], p306) and
Rukin's work [23].

\indent \indent For each $a$, let \begin{eqnarray}
\pmb{B}^{+}_{a}(\mbS_{n})
&=& \mbS^{-1}_{n}
- \left[
\begin{array}{cc}
{\bf 0}& {\bf 0}\\
{\bf 0}& \mbS_{na'a'}^{-1}
\end{array}
\right] ,
\end{eqnarray}
where $\pmb{B}^{+}_{a}(\mbS_{n})$ is positive semi-definite (p.s.d.) of rank $|a|$, ${\mbox{\reg}}$. Then we may
note that $\{(\mXbar_{n}, {\bf S}_{n})|~ n{\mXbar}'_{na:a'}{\bf S}^{-1}_{naa:a'}{\mXbar}_{na:a'}\leq k\}=\{(\mXbar_{n},
{\bf S}_{n})|~n\mXbar'_{n}B^{+}_{a}({\bf S}_{n})\mXbar_{n} \leq k\}$, for a suitable $k$, ${\mbox{\reg}}$. Tsai [29] used
the notions of the Moore-Penrose generalized inverse of matrices, matrix-convex (matrix-concave) functions, and some
related lemmas to show the following.

\vspace{0.3cm}
\noindent {\bf Proposition 4.1.} (Tsai [29]). For each $a$, the set $\{ n{\mXbar}'_{na:a'}{\bf S}^{-1}_{naa:a'}
{\mXbar}_{na:a'}\leq k\}$, for a suitable $k$, is convex in the space $({\mXbar}_{n}, {\bf S}_{n})$.
\vspace{0.2cm}

\vspace{0.3cm}
\noindent {\bf Theorem 4.5}. {\it Let ${\mathcal M}$ be the set of all the $p\times p$ positive definite matrices,
then ${\cal A}_{U}$ is closed convex on $R^{p}\times {\cal M}$.}

\noindent {\bf Proof}. Recall that we assume that ${\bf X}_{i}, i=1, \ldots, n,$ are independently multinormal distributed,
then it is easy to see that ${\mXbar}_{n}$ and ${\bf S}_{n}$ are independent, and ${\mXbar}_{na:a'}$ and ${\mXbar}_{na'}$
are independent. For each $a$, note that $I_{na}=\{({\mXbar}_{n}, {\bf S}_{n})|~{\mXbar}_{na:a'} > {\bf 0},
{\bf S}_{na'a'}^{-1}{\mXbar}_{na'} \leq {\bf 0}\}, {\mbox{\reg}}$, and $\cup_{\emptyset}^{P} I_{na}=R^{p}\times {\cal M}.$
Since that $\{(\mXbar_{n}, {\bf S}_{n})|~n{\mXbar}'_{na:a'}{\bf S}^{-1}_{naa:a'}{\mXbar}_{na:a'}\leq k\}I_{na}={\cal A}_{Ua}$, defined in (4.7),
thus by Proposition 4.1 then we have that ${\cal A}_{Ua}$ is closed convex in the space $I_{na}$, ${\mbox{\reg}}$.
Note that, any two sets ${\cal A}_{Ua}$ and ${\cal A}_{Ub}$,with $||a|-|b||=1$, are not disjoint and the intersection
of these two sets is identical to their common extreme set (Rockafellar [20]), $\emptyset \subseteq a,b \subseteq P$.
Treating ${\cal V}_{n}^{-}={A}_{U\emptyset}=\{({\mXbar}_{n}, {\bf S}_{n})|~{\bf S}_{n}^{-1}{\mXbar}_{n}\leq {\bf 0}\}
(=I_{n\emptyset})$, as the skeleton (pivotal set), we may note that ${\cal A}_{U}={\cal V}_{n}^{-}\cup_{attach}{\cal A}_{Ua}$,
where $\cup_{attach}$ means that for each $a$, $\emptyset \subset a \subseteq P$. First, the hyperspace ${\cal A}_{Ua}$
is attached to the common extreme set with the set of ${\cal V}_{n}^{-}$ from $a=\emptyset$ to $|a|=1$, and then the
hyperspace ${\cal A}_{Ub}$ (with $||b|-|a||=1$) is attached to the common extreme set of ${\cal A}_{Ua}$. Continue the
processes until $a=P$ to attach the set ${\cal A}_{UP}$. Then we may note that the set ${\cal A}_{{U}}$ is the convex
hull of ${\cal A}_{Ua}, \emptyset \subseteq a \subseteq P$. Thus the set ${\cal A}_{{U}}$ is closed convex, the probability
of $\partial {\cal A}_{{U}}$ is zero.
\vspace{0.2cm}

\vspace{0.3cm}
\noindent {\bf Theorem 4.6}. {\it Let ${\mathcal M}$ be the set of all the $p\times p$ positive definite matrices, then
${\cal A}_{U^{*}}$ is closed convex on $R^{p}\times {\cal M}$.}

\noindent {\bf Proof}. For each $a, P_{1}\subseteq a \subseteq P$, the set $\{ n{\mXbar}'_{na:a'}{\bf S}^{-1}_{naa:a'}
{\mXbar}_{na:a'} \leq k_{1} \} I_{na}$, for a suitable $k_{1}$, is closed convex in $({\mXbar}_{n},{\bf S}_{n})$. Thus the
sets ${\cal A}_{U^{*}P}$ and ${\cal A}_{U^{*}P_{1}}$ are closed convex. Notice that the set ${\cal A}_{U^{*}P}$ is a
half-ball in the space $\{({\mXbar}_{n}, {\bf S}_{n})|~{\Bar X}_{np} \geq  0\}$ and the set ${\cal A}_{U^{*}P_{1}}$ is an
unbounded closed cylinder-type convex cone in the space $\{({\mXbar}_{n}, {\bf S}_{n})|~{\Bar X}_{np} \leq  0\}$. Also the
extreme set of ${\cal A}_{U^{*}P}$ on the hyperplane $\{({\mXbar}_{n}, {\bf S}_{n})|~{\Bar {X}_{np}}=0\}$ is identical to
the one of ${\cal A}_{U^{*}P_{1}}$ on the hyperplane $\{({\mXbar}_{n}, {\bf S}_{n})|~{\Bar {X}_{np}}=0\}$. The sets
${\cal A}_{U^{*}P}$ and ${\cal A}_{U^{*}P_{1}}$ are not disjoint since ${\cal A}_{U^{*}P_{1}}={\cal A}_{U^{*}P}$ on the
hyperplane $\{({\mXbar}_{n}, {\bf S}_{n})|~{\Bar {X}_{np}}=0\}$, actually the intersection of these two sets is identical
to the extreme set of ${\cal A}_{U^{*}P}$ on the hyperplane $\{({\mXbar}_{n}, {\bf S}_{n})|~{\Bar {X}_{np}}=0\}$. The set
${\cal A}_{{U}^{*}}$ is the convex hull of ${\cal A}_{U^{*}P_{1}}$ and ${\cal A}_{U^{*}P}$. Thus the set ${\cal A}_{{U}^{*}}$
is closed convex.

\vspace{0.3cm}
\noindent {\it 4.2.1 $D-$admissibility}
\vspace{0.2cm}

\indent For the problem of testing against the positive orthant alternative, by Theorem 4.3, the accept region
${\cal A}_{U}$ of UIT is closed convex on $({\mXbar}_{n}, {\bf S}_{n})$. Hence we may assume that ${\cal A}_{U}$ is
disjoint with the following half-space
\begin{align}
 {\mathcal H}_{c} = \{({\mXbar}_{n}, {\bf S}_{n}) |~ {\mbomega}^{'}{\bf y}
   =n{\mbnu}^{'}{\mXbar}_{n}-\frac{1}{2}tr{\mbSigma}^{-1}{\bf S}_{n} > c \}.
\end{align}
The problem of testing against the half-space alternative can be handled similarly, we omit the details. Consider the
subspace
\begin{eqnarray}
 {\mathcal H}^{*}_{c} = \{({\mXbar}_{n}, {\bf S}_{n}) |~ {\mbomega}^{'}{\bf y}
    =n{\mbnu}^{'}{\mXbar}_{n}-\frac{1}{2} tr{\mbSigma}^{-1}({\bf S}_{n}+n{\mXbar}_{n}{\mXbar}^{'}_{n}) > c \}.
\end{eqnarray}
Note that ${\bf S}_{n}$ is p.d. with probability one,
thus $(\mXbar, {\bf S}_{n}) \in {\mathcal H}^{*}_{c}$ implies that $({\mXbar}_{n}, {\bf S}_{n}) \in {\mathcal H}_{c}$.
Hence, we have
\begin{eqnarray}
{\mathcal H}^{*}_{c} \subseteq {\mathcal H}_{c}.
\end{eqnarray}
Consequently, the intersection of ${\mathcal A}_{{U}}$ and ${\mathcal H}^{*}_{c}$ is empty too.

\indent Our aim is to show that ${\mbnu} \in {\cal O}^{+}_{p}\backslash \{{\bf 0}\}$. Suppose that the statement is not true,
then we have to show that ${\mbnu}={\bf 0}$ leads to a contradiction.

\indent Take ${\bf S}_{n}={\mbSigma}$ and ${\mbnu}={\bf 0}$, then we have
\begin{eqnarray}
{\mathcal H}^{*}_{c}=\{({\mXbar}_{n}, {\bf S}_{n}) \mid~n{\mXbar}_{n}^{'}{\mbS}_{n}^{-1}{\mXbar}_{n} \leq -(p+2c)\}.
\end{eqnarray}

\indent In this case, one must have $p+2c < 0$ to ensure that ${\mathcal H}^{*}_{c}$ is not an empty set. Also note that
given ${\bf S}_{n}={\mbSigma}$, ${\cal A}_{T^{2}, t^{2}_{\alpha}}$ becomes
\begin{eqnarray}
{\cal A}_{T^{2}, t^{2}_{\alpha}}=\{({\mXbar}_{n}, {\bf S}_{n}) \mid~n{\mXbar}_{n}^{'}{\mbS}_{n}^{-1}{\mXbar}_{n} \leq
t^{2}_{\alpha} \}.
\end{eqnarray}

\indent It is obvious that ${\mathcal H}^{*}_{c} \cap {\cal A}_{T^{2}, t^{2}_{\alpha}} \ne \emptyset $, and hence
${\mathcal H}^{*}_{c} \cap {\cal A}_{U} \ne \emptyset $. This contradicts the fact that (4.5) and (4.12)
are disjoint.

\indent Next, we study whether the condition (ii) of Stein Theorem (see Theorem 5.6.5 of Anderson [1]) holds or not.
For each fixed $a, \emptyset \subseteq a \subseteq P$, let
\begin{align}
\mbg = \left[
\begin{array}{cc}
{\bf I}& -\mbSigma_{aa^{'}}\mbSigma_{a^{'}a^{'}}^{-1}\\
{\bf 0}& {\bf I}
\end{array}
\right], \quad
(\mbnu, \mbSigma) \stackrel{\mbg}{\longrightarrow}
(\mbg\mbnu, \mbg\mbSigma\mbg')
\end{align}
and write
\begin{align}
\widetilde{\mbSigma} &= \mbg \mbSigma \mbg'   \\  \nonumber
&= \left[
\begin{array}{cc}
\mbSigma_{aa} - \mbSigma_{aa^{'}} \mbSigma_{a^{'}a^{'}}^{-1} \mbSigma_{a^{'}a}& {\bf 0}\\
{\bf 0}& \mbSigma_{a^{'}a^{'}}
\end{array}
\right]  \\  \nonumber
&= \left[
\begin{array}{cc}
\mbSigma_{aa:a^{'}}& {\bf 0}\\
{\bf 0}& \mbSigma_{a^{'}a^{'}}
\end{array}
\right]_.
\end{align}
Notably, $\widetilde{\mbSigma} = \widetilde{\mbSigma}'$.
Let $\mZbar_n = \mbg \mXbar_n$ and $\mbS_{n0}=\mbg \mbS_n \mbg'$.
Then,
\begin{align}
\mZbar_n = \left(
\begin{array}{c}
\mXbar_{na} - \mbSigma_{aa^{'}}\mbSigma_{a^{'}a^{'}}^{-1}\mXbar_{na^{'}}\\
\mXbar_{na^{'}}
\end{array}
\right)_{,}
\end{align}
and
\begin{align}
& \tr \mbSigma^{-1}(\mbS_n + n\mXbar_n\mXbar_n')   \\  \nonumber
&= \tr (\mbg \mbSigma \mbg')^{-1}\mbg
(\mbS_n + n\mXbar_n\mXbar_n')\mbg'  \\  \nonumber
&= \tr \widetilde{\mbSigma}^{-1}(\mbS_{n0} + n\mZbar_n \mZbar'_n)_.
\end{align}
We may choose $\widetilde{\mbnu}=({\mbg}^{-1})^{'}{\mbnu}$ so that ${\mbnu}^{'}\mXbar_n=\widetilde{\mbnu}'\mZbar_n$,
where
\begin{align}
 \hspace{1cm} \widetilde{\mbnu}
&  =({\mbg}^{-1})^
{'} \mbnu      \\  \nonumber
& = \left(
\begin{array}{c}
{\mbSigma}^{-1}_{aa:a^{'}}({\mbnu}_{a}-{\mbSigma}_{aa^{'}}{\mbSigma}^{-1}_{a^{'}a^{'}}{\mbnu}_{a^{'}})\\
{\mbSigma}^{-1}_{a^{'}a^{'}}{\mbnu}_{a^{'}}
\end{array}
\right)   \\  \nonumber
&\triangleq \left(
\begin{array}{c}
\widetilde{\mbnu}_{a}\\
\widetilde{\mbnu}_{a^{'}}
\end{array}
\right).
\end{align}
Thus, the subspace ${\mathcal H}^{*a}_{c}$ becomes
\begin{align}
{\mathcal H}^{*a}_{c} & = \{\, (\mZbar_{n}, \mbS_{n0}) \mid n\widetilde{\mbnu}'\mZbar_{n}
- \frac{1}{2} \tr \widetilde{\mbSigma}^{-1}(\mbS_{n0} + n\mZbar_{n}\mZbar_{n}') > c \} \\ \nonumber
   & = \{\, (\mZbar_{n}, \mbS_{n0}) \mid n\widetilde{\mbnu}'\mZbar_{n}
 - \frac{n}{2} \mZbar_{n}^{'}\widetilde{\mbSigma}^{-1}\mZbar_{n}-\frac{1}{2}
   \tr \widetilde{\mbSigma}^{-1}\mbS_{n0}  > c \} \\ \nonumber
&  = \{\, (\mXbar_{n}, \mbS_{n}) \mid n{\mbnu}'\mXbar_{n}
   - \frac{n}{2} \mXbar^{'}_{0na:a^{'}}{\mbSigma}^{-1}_{aa:a^{'}}\mXbar_{ona:a^{'}}
  -\frac{n}{2}\mXbar_{na^{'}}^{'}{\mbSigma}^{-1}_{a^{'}a^{'}}\mXbar_{na^{'}}-\frac{1}{2}
    \tr {\mbSigma}^{-1}\mbS_{n}  > c \},
\end{align}
where $\mXbar_{ona:a^{'}}=\mXbar_{na}-{\mbSigma}_{aa^{'}}{\mbSigma}^{-1}_{a^{'}a^{'}}\mXbar_{na^{'}}$.

\indent Since ${\bf S}_{n}$ converges to ${\mbSigma}$ almost surely, thus with probability one for each fixed $a$ we have
\begin{align}
{\mathcal H}^{*a}_{c} &  = \{\, (\mXbar_{n}, \mbS_{n}) \mid n{\mbnu}'\mXbar_{n}
   - \frac{n}{2} \mXbar^{'}_{na:a^{'}}{\mbS}^{-1}_{naa:a^{'}}\mXbar_{na:a^{'}}
  -\frac{n}{2}\mXbar_{na^{'}}^{'}{\mbS}^{-1}_{na^{'}a^{'}}\mXbar_{na^{'}} > c+\frac{p}{2}\}, \\ \nonumber
  &=\{\, (\mXbar_{n}, \mbS_{n}) \mid n{\mbnu}^{'}_{a}\mXbar_{na:a^{'}}
   - \frac{n}{2} \mXbar^{'}_{na:a^{'}}{\mbS}^{-1}_{naa:a^{'}}\mXbar_{na:a^{'}}
  +\frac{n}{2}[\mbnu^{'}_{a^{'}}{\mbSigma}_{a^{'}a^{'}}\mbnu_{a^{'}}-(\mXbar_{na^{'}}  \\ \nonumber
  & ~~~~~~-{\mbSigma}_{a^{'}a^{'}}\mbnu_{a^{'}})^{'}{\mbSigma}^{-1}_{a^{'}a^{'}}
   (\mXbar_{na^{'}}-{\mbSigma}_{a^{'}a^{'}}\mbnu_{a^{'}})] > c+\frac{p}{2}\},  \nonumber
\end{align}

\indent Take $\mbnu_{a^{'}}={\bf 0}$ and ${\mbSigma}={\mbS}_{n}$, then for each fixed $a, \emptyset \subseteq a
\subseteq P$, we have
\begin{align}
{\mathcal H}^{*a}_{c} &=\{\, (\mXbar_{n}, \mbS_{n}) \mid n{\mbnu}^{'}_{a}\mXbar_{na:a^{'}}- \frac{n}{2}
   \mXbar^{'}_{na:a^{'}}{\mbS}^{-1}_{naa:a^{'}}\mXbar_{na:a^{'}}-\frac{n}{2}\mXbar^{'}_{na^{'}}
  {\mbS}^{-1}_{na^{'}a^{'}} \mXbar_{na^{'}} > c+\frac{p}{2}\}.
\end{align}

\indent By Theorem 4.3, we assume that ${\mathcal H}_{c} \cap {\cal A}_{U}=\emptyset$ in (4.13). Note that ${\cal A}_{{U}}
=\cup_{\emptyset}^{P}{\cal A}_{Ua}$ which is defined in (4.6) and (4.7). For each fixed $a$, we may notice that
${\mathcal H}^{*a}_{c} \cap {\cal A}_{Ua}=\emptyset, \emptyset \subseteq a \subseteq P$. Thus, for each fixed $a$, there
exist $\mbnu_{a} \geq {\bf 0}, \mbnu_{a} \ne {\bf 0}$, and $\mbnu_{a^{'}}={\bf 0}, \emptyset \subseteq a \subseteq P$.
Let ${\cal O}_{a}=\{\mbnu \in R^{p} \mid \mbnu_{a}\geq {\bf 0}, {\mbnu}_{a}\ne {\bf 0}$,  and $\mbnu_{a^{'}}={\bf 0} )\}$,
then we have ${\cal O}^{+}_{p}=\cup_{\emptyset}^{P}{\cal O}_{a}$. It means that there exists $\mbomega_{1} \in
{\cal O}^{+}_{p}$ such that for arbitrarily large $\lambda$ the vector ${\mbomega}_{1}+\lambda {\mbomega} \in
{\cal O}^{+}_{p} \backslash \{{\bf 0}\}$. Namely, the condition (ii) of Theorem 5.6.5 of Anderson [1] holds. Hence by Theorem
8 and Corollary 2 of Lehmann ([13], p306-308), we may conclude that for the problem (1.1) both the UITs are $d$-admissible.

\vspace{0.3cm}
\noindent {\bf Theorem 4.7}. {\it For the problem of testing $H_{0}: {\mbtheta} = {\bf 0}$
against $H_{1 {\cal O}^{+}}: {\mbtheta} \in {\cal O}_{p}^{+} \backslash \{{\bf 0}\}$, the UIT is $d$-admissible.}
\vspace{0.2cm}

\indent As for the problem of testing against the half-space alternative, the $d$-admissibility of the UIT can be similarly
proved the same as that of Theorem 4.7 by only considering the cases $a=P_{1}$ and $a=P$. Hence, we omit the details.

\vspace{0.3cm}
\noindent {\it 4.2.2 $\alpha-$admissibility}
\vspace{0.2cm}

\indent Furthermore, with similar arguments as in Theorems 2.1 and 3.1 of Sen and Tsai [24], then we may also show that for
the problem of testing against the half-space alternative, the UIT is similar and unbiased. Therefore, by Corollary 2 of
Lehmann [13] we have the following.

\vspace{0.3cm}
\noindent {\bf Theorem 4.8}. {\it For the problem of testing $H_{0}: {\mbtheta} = {\bf 0}$ against $H_{1 {\cal H}^{*}}:
{\mbtheta} \in {\cal H}_{p}^{*} \backslash \{{\bf 0}\}$, the UIT is $\alpha$-admissible.}
\vspace{0.2cm}

\vspace{0.3cm}
\def \theequation{5.\arabic{equation}}
\setcounter{equation}{0}

\noindent {\bf 5. The general remarks}
\vspace{0.2cm}

\indent Fisher's approach to reporting $p$-values involves identifying the optimal method of the inside information under
the null hypothesis. Neyman-Pearson's theory aims to identify the optimal power properties under the alternative hypothesis. Both the
Fisher's approach and Neyman-Pearson's theory are compelling for addressing the problems of testing against both the
unrestricted alternative and the half-space alternative, while it is not pertain to the problem of testing against the
positive orthant space alternative world. For the issue of testing against the positive orthant space alternative, both
the null distributions of LRT and UIT depend on the unknown covariance matrix; as such, we have difficulty adopting the
Fisher's method to reporting $p$-values. Meanwhile, the Neyman-Pearson's optimal power theory indicates that both the LRT
and the UIT are $\alpha$-inadmissible as proved in Propositions 2.1 and 2.2. For the problem (1.1) of testing against
restricted alternatives, Perlman and Wu [19] downplayed the significance of theoretical properties, such as the
$\alpha$-admissibility, similarity and unbiasedness of Neyman-Pearson theory, in favor of the LRT. Berger [5] criticized
Perlman and Wu's suggestion lacks any substantial support for the LRT. One of the preliminary objective is to seek a
potential resolution for the arguments between Perlman and Wu [19] and Berger [5].

\indent As we know from Section 2 the power domination problems of hypothesis testing heavily depend on the choice of
the critical points that are used to report the $p$-values. Sometimes, the danger stems from too narrow definitions
of what is meant by optimality and the strongest intuition can sometimes go astray. Hence, the question naturally
raised is whether there exists any satisfactory, unified approach to overcome the difficulties. The spirit of compromise
between Fisher's approach and Neyman-Pearson's optimal theory without detailed consideration of power may shed light on the
testing hypothesis theory. Imposing the balance between type 1 error and power, Wald's decision theory may pave a unified
way to combine the best features of Neyman-Pearson's and Fisher's ideas. Hence, we suggest that power comparisons of
tests would make sense only because all the compared tests are similar. As such, Wald's decision theory
via $d$-admissibility and Neyman-Pearson's optimal theory is essentially one theory, not two. For the problem (1.1), we
show that UITs are $d$-admissible. Thus, a framework to achieve the goal of reconciliation of Bayesian, frequentist
(Neyman-Pearson-Wald approach) and Fisherian paradigms for the problem of testing against restricted alternatives
(closed convex cones) can then be established. With this new novel framework, to the problem of testing against the positive
orthant space alternative both the theoretical properties of $\alpha-$inadmissibility, similarity and unbiasedness of
Neyman-Pearson theory and the method of Fisher's reporting $p$-values will be simultaneously downplayed to get rid of
those phenomena against common statistical sense. This new framework might provide a chance to develop a new theory of
hypothesis testing problems.

\indent Recall that ${\mathcal A}_{L}=\cup_{\emptyset}^{P}{\cal A}_{La}$ is the acceptance region of the LRT for the problem
of testing against the positive orthant space alternative, i.e., $H_{0}: {\mbtheta} = {\bf 0}$ against $H_{1}: {\mbtheta}
\in {\cal O}^{+}_{p} \backslash \{{\bf 0}\}$, where ${\cal A}_{La}=\{({\mXbar}_{n}, {\bf S}_{n})|~{\frac{n{\mXbar}'_{na:a'}
{\bf S}_{naa:a'}^{-1} {\mXbar}_{na:a'}}{1+n{\mXbar}'_{na'} {\bf S}_{na'a'}^ {-1} {\mXbar}_{na'}}} \leq l_{\alpha} \}I_{na}$,
{\reg}, with $l_{\alpha}$ being the critical point of the level of significance $\alpha$. For testing against the global
alternative, the LRT always plays an important role in the history of statistical hypothesis testing. As we have observed in the discussions, the LRT can be considered $\alpha$-inadmissible when testing against the
positive orthant space alternative, the LRT is $\alpha$-inadmissible. So does the $\alpha$-inadmissible result for the UIT.
Traditionally, the Birnbaum-Stein method stipulates convexity for the acceptance regions of tests. In this paper, we
generalize the Birnbaum-Stein method to show the acceptance regions of UITs are convex, and further to demonstrate that the
UITs are $d$-admissible for the problem (1.1). It is important to note that when testing against the positive orthant space
alternative, the acceptance region ${\mathcal A}_{L}$ of the LRT is no longer a convex set. Similarly, for testing against
the half-space alternative, the acceptance region of the LRT is not a convex set. We encounter with the situation that the
acceptance regions of LRTs are hyperbolic sets rather than the convex sets for the problem (1.1), this motivates future
studies to investigate whether the LRTs are $d$-admissible for the problem.

\indent In Section 4, a framework for UIT with Wald's $d$-admissibility is established for the problem of testing against
the positive orthant space alternative, assuming the Birnbaum-Stein method, where the acceptance region is a closed
convex set. It's crucial to note that the acceptance region of the LRT for testing against of the positive orthant space
alternative is non-convex, as unequivocally discussed in the previous paragraph. Traditionally, the LRT is important for
testing against the unrestricted alternative. However, it has been irrefutably shown that both the LRT and the UIT are
biased and $\alpha$-inadmissible for testing against the positive orthant space alternative. Despite this, the UIT is
unequivocally shown to be $d$-admissible, downplaying the effects of $\alpha$-inadmissibility and bias. These phenomena
are due to the null distribution function depending on the unknown covariance matrix. It's unquestionably reasonable to
expect that the LRT is also $d$-admissible, just like the UIT, for testing against the positive orthant space alternative.
This raises the question whether the assumption that the acceptance region is a closed convex set can be relaxed so
that the LRT is $d$-admissible. It is vital to explore the possibility for further establishing a framework for the
LRT for the problem of testing against the positive orthant space alternative.

\vspace{0.3cm}
\vspace{0.2cm}

\vspace {0.3cm}
\noindent {\bf References}
\vspace{0.2cm}

\begin{enumerate}
\item  {T.W. Anderson},  An Introduction to Multivariate Statistical Analysis. 2nd edition. Wiely, New York, 1984.

\item  {D. Basu}, On the elimination of nuisance parameters.  {J. Amer. Statist. Assoc.}  {72} (1977) 355-366.

\item  {J.O. Berger, B. Liseo, R.L. Wolpert}, Integrated likelihood methods for eliminating nuisance parameters.
                {Statist. Sci.} {14} (1999) 1-28.

\item  {R.L. Berger},  Uniformly more powerful tests for hypotheses concerning linear inequalities and normal means.
              {J. Amer. Statist. Assoc.} {84} (1989) 192-199.

\item  {R.L. Berger}, Comment on ``The emperor's new tests". {Statist. Sci.} {14} (1999) 370-373.

\item  {A. Birnbaum}, Characterization of complete classes of tests of some multiparametric hypothesis, with
             applications to likelihood ratio tests. {Ann. Math. Statist.} {26} (1955) 21-36.

\item  {L.D. Brown}, {Fundamentals of Statistical Exponential families, with Application in Statistical Decision Theory},
           Institute of Mathematical Statistics, Lecture Notes-Monograph Series, Volume {9}. Hayward, California. 1986.

\item  {T. Chang, D. Eaves}, Reference priors for the orbit in a group model. {Ann. Statist.} {18} (1990) 1595-1614.

\item  {M.L. Eaton}, A complete class theorem for multidimensional one sided alternatives. {Ann. Math. Statist.} {41}
             (1970) 1884-1888.

\item  {B. Efron}, R.A. Fisher in the 21st century, Statist. Sci. {13} (1998) 95-144.

\item  {J. Kiefer, R. Schwartz}, Admissible Bayes character of $T^2$-, $R^2$-, and other fully invariant tests
           for classical multivariate normal problems", {Ann. Math. Statist.} {36} (1965) 747-770.

\item  {E.L. Lehmann}, The Fisher, Neyman-Pearson theories of testing hypotheses: One theory or two? {J. Amer. Statist.
           Assoc.} {88} (1993) 1242-1249.

\item  {E.L. Lehmann}, {Testing Statistical Hypotheses}. 2nd edition. New York: Wiley, 1986.

\item  {J.I. Marden},  Minimal complete classes of tests of hypothesis with multivariate one-sided alternatives.
          {Ann. Statist.} {10} (1982) 962-970.

\item  {J.A. Menendez, B. Salvador}, Anomalies of the likelihood ratio test for testing restricted hypotheses.
              {Ann. Statist.} {19} (1991) 889-898.

\item  {J.A. Menendez, C. Rueda, B. Salvador}, Dominance of likelihood ratio tests under cone constraints.
                {Ann. Statist.} {20} (1992) 2087-2099.

\item  {P.E. N\"{u}esch}, On the problem of testing location in multivariate problems for restricted alternatives.
             {Ann. Math. Statist.} {37} (1966) 113-119.

\item  {M.D. Perlman}, One-sided problems in multivariate analysis. {Ann. Math. Statist.} {40} (1969) 549-567.

\item  {M.D. Perlman, L. Wu}, The emperor's new tests. {Statist. Sci.} {14} (1999) 355-369.

\item  {R.T. Rockafellar}, {Convex Analysis}. Princeton University Press, 1972.

\item  {S.N. Roy}, On a heuristic method of test construction and its use in multivariate analysis.
              {Ann. Math. Statist.} {24} (1953) 220-238.

\item  {S.N. Roy, R. Gnanadesikan, J.N. Srivastava},  {Analysis and Deign of Certain Quantitative Multiresponse Experiments}.
             Pergamon Press, Oxford, 1972

\item  {A.L. Rukhin}, Admissibility: Survey of a concept in progress. Int. Statist. Rev. {63} (1995) 95-115.

\item  {P.K. Sen, M.-T. Tsai}, (1999), Two-stage LR and UI tests for one-sided alternatives: multivariate mean with
             nuisance dispersion matrix. {J. Multivariate Anal.} {68} (1999) 264-282.

\item  {J.B. Simaika}, On an optimum property of two important statistical tests. {Biometrika}  {32}, (1941) 70-80.

\item  {C. Stein}, A two-sample test for a linear hypothesis whose power function is independent of $\sigma$.
              {Ann. Math. Statist.} {16} (1945) 243-258.

\item  {C. Stein}, The admissibility of Hotelling's $T^{2}$-test. {Ann. Math. Statist.} {27} (1956) 616-623.

\item  {D.I. Tang},  Uniformly more powerful tests in a one-sided multivariate problem. {J. Amer. Statist. Assoc.}
              {89} (1994) 1006-1011.

\item  {M.-T. Tsai}, Admissibility of invariant tests for means with covariates. {Math. Method. Statist.}
          {28} (2019) 243-261.

\item  {M.-T. Tsai, P.K. Sen}, On inadmissibility of Hotelling's $T^{2}$-test for restricted alternatives.
           {J. Multivariate Anal.} {89} (2004) 87-96.

\item  {A. Wald}, {Statistical Decision Functions}. Wiely, New York, 1950.

\item   {Y. Wang, M.P. McDermott}, Conditional likelihood ratio test for a nonnegative normal mean vector.
             {J. Amer. Statist. Assoc.} {93} (1998) 380-386.

\end{enumerate}

\vspace{0.3cm}
\newlength{\lewidth}
\newlength{\riwidth}
\settowidth{\lewidth}{University of Virginia at Charlottesville, Virginia}
\settowidth{\riwidth}{Institute of Statistical Science at Taipei}


\end{document}